\documentclass[10pt,twoside]{article}
\usepackage[latin1]{inputenc}
\usepackage{amsmath}
\usepackage{graphicx}
\usepackage{yhmath}
\usepackage[table]{xcolor}
\usepackage{ mathrsfs} 
\usepackage{amssymb}
\usepackage{url}

\input epsf
\def\limiten{\renewcommand{\arraystretch}{0.5}
\begin{array}[t]{c}\stackrel{}{\longrightarrow} \\
{\scriptstyle n\rightarrow
\infty}\end{array}\renewcommand{\arraystretch}{1}}

\def\limitepsn{\renewcommand{\arraystretch}{0.5}
\begin{array}[t]{c}\stackrel{a.s.}{\longrightarrow} \\
{\scriptstyle n \rightarrow
\infty}\end{array}\renewcommand{\arraystretch}{1}}

\def\limitepsm{\renewcommand{\arraystretch}{0.5}
\begin{array}[t]{c}\stackrel{a.s.}{\longrightarrow} \\
{\scriptstyle m \rightarrow
\infty}\end{array}\renewcommand{\arraystretch}{1}}

\DeclareMathOperator\argmax{argmax}

\renewcommand{\Box}{\hfill\rule{0.25cm}{0.25cm}} 

\newcommand\esp[1]{\,\,\, \text{#1}\,\,\,}

\newtheorem{Prop}{Proposition}[section]
\newtheorem{lem}{Lemma}[section]
\newtheorem{Theo}{Theorem}[section]
\newtheorem{cor}{Corollary}[section]

\newenvironment{dem}{\ \\ {\bf Proof. }}
{\Box\par\medskip\noindent}

\newcommand{\E}{\ensuremath{\mathbb{E}}}
\newcommand{\R}{\ensuremath{\mathbb{R}}}
\newcommand{\Z}{\ensuremath{\mathbb{Z}}}

\newcommand{\N}{\ensuremath{\mathbb{N}}}

\definecolor{grisclair}{gray}{0.9}

\renewcommand{\arraystretch}{.8}

\begin{document}
\title{\bf Inference for nonstationary time series of counts with application to change-point problems}
 \maketitle \vspace{-1.0cm}
 \begin{center}
   William Kengne$^{\text{a}}$ \footnote{Developed within the ANR BREAKRISK (ANR-17-CE26-0001-01) } and
  Isidore Séraphin Ngongo$^{\text{b}}$. 
 \end{center}

  \begin{center}
 $^{\text{a}}$ CY Cergy Paris Université, CNRS, THEMA, F-95000 Cergy, France.\\
{\it $^{\text{b}}$  ENS, Université de Yaoundé 1, Cameroun.\\
  E-mail:  william.kengne@u-cergy.fr  ;  ngongoisidore@gmail.com}
 \end{center}

 \pagestyle{myheadings}
 \markboth{Nonstationary time series of counts with application to change-point}{W. Kengne and I. Ngongo}

\textbf{Abstract} :
  We consider an integer-valued time series $Y=(Y_t)_{t\in\Z}$ where the models after a time $k^*$ is Poisson autoregressive with the conditional mean that depends
  on a parameter $\theta^*\in\Theta\subset\R^d$. The structure of the process before $k^*$ is unknown; 
  it could be any other integer-valued time series, that is, the process $Y$ could be nonstationary.
  It is established that the maximum likelihood estimator of $\theta^*$ computed on the nonstationary observations is consistent and asymptotically normal.
  Next, we carry out  the sequential change-point detection in a large class of Poisson autoregressive models.  
  We propose a monitoring scheme for detecting change in the model.
  The procedure is based on an updated estimator which is computed without the historical observations.
 The asymptotic behavior of the detector is studied, in particular, the above result on the inference in a nonstationary setting are applied to prove that the proposed procedure
   is consistent.
 A simulation study as well as a real data application  are provided.
  \\


\section{Introduction}\label{intro}
We consider a process $Y=(Y_t)_{t \in \Z}$ satisfying
\begin{equation}\label{NonStmodel}
   Y_t | \mathcal{F}_{t-1}   \sim\mbox{ Poisson}(\lambda_t)  ~ \text{with} ~ \lambda_t = \E(Y_t | \mathcal{F}_{t-1});
   \end{equation}
 where $\mathcal{F}_t=\sigma(Y_s, ~ s \leq t)$ is the $\sigma$-field generated by the whole past.
 A large literature on this model has recently been developed by assuming that $\lambda_t = \E(Y_t | \mathcal{F}_{t-1}) = f(Y_{t-1}, Y_{t-2},\cdots)$ for all $t \in \Z$,
 where $f$ is a measurable non-negative function, satisfying some Lipschitz-type conditions.
 This entails that the process $(Y_t , \lambda_t)_{t \in \Z}$ is strict stationary  with finite moment of any order.
 But, such result does not hold in many practical situations. For instance, in the change-point problem, it often occurs that
   \begin{equation*}
   \lambda_t = \begin{cases}
                f_0(Y_{t-1}, Y_{t-2},\cdots) & \text{for } t \leq   k^* \\
                f_1(Y_{t-1}, Y_{t-2},\cdots)  &  \text{for } t  > k^*
  \end{cases}
   \end{equation*}
 with $f_0 \neq f_1$ and $k^* \in \Z$. Thus, the process $(Y_t , \lambda_t)_{t \in \Z}$ is not stationary.

 \medskip

  We consider a nonstationary autoregressive process $Y=(Y_t)_{t \in \Z}$ in a parametric framework; we assume that $Y$ satisfying
  \begin{equation}\label{NonStmodel_par}
      Y_t | \mathcal{F}_{t-1}   \sim\mbox{ Poisson}(\lambda_t)  ~ \text{with} ~   \lambda_t = \E(Y_t | \mathcal{F}_{t-1}) = f_{\theta^*}(Y_{t-1}, Y_{t-2},\cdots)    \text{ for all } t >  k^*;
   \end{equation}
 with $k^* \in \Z$, $\theta^*$ is the parameter belonging to a compact set $\Theta \subset \R^d$ ($d \in \N$) and $f_{\theta}$ a measurable non-negative function,
   assumed to be known up to the parameter $\theta$. If (\ref{NonStmodel_par}) holds for $t < k^*$, then with some Lipschitz-type conditions on $f$, the process
  $(Y_t , \lambda_t)_{t \in \Z}$ is strict stationary  with finite moment of any order (see for instance Doukhan {\it et al.} (2012)).
  We focus here in a more general situation where the structure of the process $(Y_t)_{t\leq k^*}$ is assumed to be unknown;
  it could be a Poisson autoregressive model depending on a parameter different from $\theta^*$ or could be any other integer-valued time series.

  \medskip

  In this work, we firstly study the inference on the parameter $\theta^*$ in the model (\ref{NonStmodel_par}). This task has been considered by several authors;
  see among others Fokianos {\it et al.} (2009), Fokianos and Tj{\o}stheim (2012), Doukhan and Kengne (2015).
  These works (and many other) have been developed under the assumption that the process  $(Y_t)_{t \in \Z}$ is strict stationary ;
  which restrict the application area of such results.
   To deal with the model (\ref{NonStmodel_par}), we conduct some preliminary works that leads to approximate the nonstationary process with its stationary regime.
   Under some classical Lipschitz-type condition on the function $f$, there exists (see \cite{Doukhan2012,Doukhan2013corr}) a strict
  stationary process $\tilde{Y} =(\tilde{Y}_t)_{t \in \Z}$  with finite moment of any order, satisfying :
   \begin{equation}\label{modYtilde}
   \tilde{Y}_t | \mathcal{\tilde{F}}_{t-1}   \sim\mbox{ Poisson}(\tilde{\lambda}_t) ~ \text{with} ~  \tilde{\lambda}_t=f_{\theta^*}(\tilde{Y}_{t-1}, \tilde{Y}_{t-2},\cdots) ~ ~
   \text{for } t \in \Z
     \end{equation}
 where $ \mathcal{\tilde{F}}_{t} = \sigma(\tilde{Y}_s, s\leq t)$ is the $\sigma$-field generated by the whole past of $\tilde{Y}$. \\
 \noindent Let us remark that, models (\ref{NonStmodel}), (\ref{NonStmodel_par}) and (\ref{modYtilde})  can be represented in terms of Poisson processes.
    Let $\{ N_t(\cdot) \esp{;} t=1,2,\cdots \}$ be a sequence of independent Poisson processes of unit intensity. $Y_t$ and $\tilde{Y}_t$ can respectively be seen as the number
    (say $N_t(\lambda_t)$) of events of $N_t(\cdot)$ that occurs in the time interval $[0,\lambda_t]$ and $[0, \tilde{\lambda}_t]$.
    Therefore, we can also write
     \begin{equation}\label{modelPP}
   Y_t = N_t(\lambda_t), ~  \tilde{Y}_t = N_t(\tilde{\lambda}_t)  \text{ with } \lambda_t = f_{\theta^*}( Y_{t-1},Y_{t-2} \ldots) \text{ and }
    \tilde{\lambda}_t=f_{\theta^*}(\tilde{Y}_{t-1}, \tilde{Y}_{t-2},\cdots) \text{ for all } t > k^*.
   \end{equation}
  This representation is useful to approximate the processes $(Y_t)_{t \geq k^*}$ and  $(\tilde{Y}_t)_{t \geq k^*}$.
  The question of this approximation has been addressed by Doukhan and Kengne \cite{Doukhan2015} (see Remark 4.1). In this work, we provide a detailed proof of this problem.
  In particular, we show that the expectation $\E| Y_{k^*+\ell} - \tilde{Y}_{k^*+\ell} |$ (for $\ell \geq 1$) can be controlled and tends to zero when $\ell$  goes to infinity, see Lemma \ref{lem2}.
  These approximation results are applied to establish that the conditional maximum likelihood estimator (MLE) of $\theta^*$, based on the nonstationary observations is consistent and asymptotically normal.
 Also, let us stress that numerous papers on change-point problem assume that the process is stationary after the breakpoint; see for instance Doukhan and Kengne \cite{Doukhan2015}, Diop and Kengne (2017), Franke {\it et al.} (2012), Kirch and Tadjuidje Kamgaing (2015). 
 This paper provides tools to avoid such condition which is quite  restrictive in practice. 

  \medskip

 As a second contribution, we consider the structural change-point problem in Poisson autoregressive models. In the retrospective (or off-line) framework, this issue has already been
 addressed. See for instance Franke {\it et al.} (2012), Kang and Lee (2014), Doukhan and Kengne \cite{Doukhan2015}, Diop and Kengne \cite{Kengne2017}.
 But these works suffer from a drawback : the (asymptotic) study under the presence on change-point is either missing or done with the stationarity assumption on the observations
 after breakpoint, which is unrealistic in many practical problems.
 For procedure proposed by these authors, stationarity assumption after change-point can been relaxed by applying Theorem \ref{theo_consistent} (see below) to get the consistency under the alternative
 of change occurs in the model. In the sequel, we focus on sequential (or on-line) framework. \\
 \noindent Assume that the process $Y=(Y_t)_{t \in \Z}$ satisfying
  \begin{equation}\label{model}
   Y_t/\mathcal{F}_{t-1}   \sim\mbox{ Poisson}(\lambda_t) ~ \text{with} ~
   \lambda_t = \begin{cases}
                f_{\theta_0^*}(Y_{t-1}, Y_{t-2},\cdots) & \text{for } t \leq  k^* \\
                f_{\theta_1^*}(Y_{t-1}, Y_{t-2},\cdots)  &  \text{for } t  > k^*
  \end{cases}
   \end{equation}
 where $\theta_0^*, \theta_1^*$ are the parameters belonging to a compact set $\Theta \subset \R^d$ ($d \in \N$) and $k^*$ a positive integer representing the possible breakpoint.
 If $\theta_0^* \neq \theta_1^*$, then a structural
   change occurs at time $ k^*$ ; otherwise, no change has occurred and the model (\ref{model}) can be simply written as
   \begin{equation}\label{model2}
   Y_t/\mathcal{F}_{t-1}   \sim\mbox{ Poisson}(\lambda_t) ~ \text{with} ~  \lambda_t=f_{\theta_0^*}(Y_{t-1}, Y_{t-2},\cdots) ~ ~  \text{for } t \in \Z .
   \end{equation}
 \indent  We follow the paradigm of Chu {\it et al.} (1996). The general idea is to use the observations
  $(Y_1,\cdots,Y_m)$ (called the historical data) that depends on the parameter $\theta_0^*$, then, one can monitor the change in the model's parameter sequentially
  from the date $m+1$ and trigger an alarm when a change is detected; by  ensuring that the probability of false alarm does not exceed a fixed level $\alpha$.
  More precisely, $k^* > m$ and $(Y_1,\cdots,Y_m)$ is assumed to be generated from the model (\ref{model}), depending on $\theta_0^*$ (without change);
  we are going to observe new data
  $Y_{m+1}, Y_{m+2},\cdots, Y_{m+k}, \cdots$. For each new observation $Y_{m+k}$, we will like to know if it is generated from a model depending on $\theta_0^*$ or from a model
  depending on $\theta_1^*$, with $\theta_0^* \neq \theta_1^*$.  This problem can be treated as a classical hypothesis testing : \\

 \noindent  $\mathbf{H_0}$: $\theta^*_0$ is constant over the observations $Y_1,\cdots,Y_m,Y_{m+1},\cdots$ {\it i.e.}  $(Y_t)_{t\in \N}$
             satisfying (\ref{model}) with  $\theta^*_0 = \theta^*_1$;\\
 ~\\
 $\mathbf{H_1}$ : the process $(Y_t)_{t\in \N}$ satisfying (\ref{model}) with  $\theta^*_0 \neq \theta^*_1$.\\

  Numerous works have been done in the sequential change-point detection according to such paradigm.
 See  among others papers,  Horváth {\it et al.} (2004)\nocite{Horvath2004}, Gombay and Serban (2009)\nocite{Gombay2009}, Na {\it et al.} (2011)\nocite{Na2011}
 Bardet and Kengne (2014)\nocite{Bardet2014} for several tests procedure for sequential change detection in a general class of time series models, including
  linear and GARCH-type models.
  Kengne (2015)\nocite{Kengne2015} proposed a fluctuation-type test procedure for sequential change detection in a large class of Poisson autoregressive model.
 %
%
%
%
 Recently,  Kirch and Tadjuidje Kamgaing \cite{Kirch2015} and Kirch and Weber (2018) have considered a large class of models (that including continuous and discrete valued time series) and
 developed a general setup based on estimating functions for sequential change-point detection.
 Estimating functions is a general estimation method and some classical procedure such as likelihood estimator, least square estimator,$\cdots$ can be treated in many cases as
 a particular class of estimating functions.
 It is well-known (Godambe (1960)\nocite{Godambe1960}) that the optimal estimating function in several classical  parametric  model is based on the score function.
 In the case of infinite memory process considered here, a more complex class of estimating functions is needed; this involves some difficulties in the application of their procedure.
 Moreover, Kirch and Tadjuidje Kamgaing \cite{Kirch2015} and Kirch and Weber \cite{Kirch2018}) impose some regularity conditions on the process after the change-point.
 These conditions, which is not easy to verify in general, are somewhere sufficient  to unify the treatment of the large class models that they have considered.

 \medskip

 We carry out a sequential test in the spirit of Bardet and Kengne \cite{Bardet2014}, and propose an open-end and closed-end (see below) procedure for monitoring changes in
 the model (\ref{model}).
 We develop a procedure where the recursive estimator is computed without the historical observations.
 It is shown that the detector converges to a well-known distribution under the null hypothesis.
 Under the alternative, we do not need any additional assumption on the process after the change-point. The consistency of the procedure is established event in the nonstationary
 setting (the previous study on inference in nonstationary models play a key role in the proof of this result).  Moreover, the test developed here is intended to early detect change than the aforementioned procedure, since it has displayed a detection delay that can be bounded by
 $\mathcal{O}_P(m^{1/2 + \epsilon})$ for any $\epsilon >0$.

 \medskip

 In the following Section 2, some classical assumptions on the model (\ref{NonStmodel_par}) as well as some examples are provided.
 The inference in the nonstationary process $Y$ is conducted in Section 3.
 Section 4 focuses on the sequential change-point detection.  Some numerical results are displayed in Section 5, whereas  Section 6 is devoted to a concluding remarks. 
 The proofs of the main results are provided in Section 7.

 \section{Assumptions and examples}
 \subsection{Assumptions}
  We will use the following classical notations:
 \begin{enumerate}
    \item $\|y\|:=\sum\limits_{j=1}^{p} |y_j|$ for any $y \in \R^p$;
    \item for any compact set $ \mathcal{K} \subseteq\R^d$ and for any function
    $g:\mathcal{K} \longrightarrow\R^{d'}$, $ \|g\|_\mathcal{K} =\sup_{\theta\in \mathcal{K} }(\|g(\theta)\|)$;
    \item for any set $ \mathcal{K} \subseteq\R^d$,  $ \ring{\mathcal{K}}$ denotes the interior of $\mathcal{K} $;
    \item $\N = \{1,2,3,\cdots  \}$ and $\N_0 = \{0,1,2,3,\cdots  \}$.
 \end{enumerate}
  \noindent Throughout the sequel, we will assume that the function $\theta \mapsto f_\theta$ is twice continuously differentiable on $\Theta$  and we need the following conditions on the model (\ref{NonStmodel_par}).

 \medskip
  \noindent For $i=0,\, 1, \, 2$ , define \\
 {\bf Assumption A$_i (\Theta) $}: $\|{\partial^i f_\theta(0)}/{\partial\theta^i}\|_\Theta<\infty$
 and there exists a sequence of non-negative real numbers $(\alpha^{(i)}_k )_{k\geq 1}$ satisfying
   $ \sum\limits_{j=1}^{\infty} \alpha^{(0)}_k  < 1$ (when $i=0$) and  $ \sum\limits_{j=1}^{\infty} \alpha^{(i)}_k <\infty$ (when $i=1,2$)
   such that
 \begin{equation*}
 \Big\|\dfrac{\partial^i f_\theta(y)}{\partial\theta^i}-\dfrac{\partial^i f_\theta(y')}{\partial\theta^i}\Big\|_{\Theta}
 \leq \sum\limits_{k=1}^{\infty}\alpha^{(i)}_k |y_k-y'_k| \quad \mbox{ for all } y, y' \in (\R^+)^{\N}.
 \end{equation*}
 Under the assumption \textbf{A}$_0 (\Theta)$,  Doukhan \textit{et al.} (2012, 2013)\nocite{Doukhan2012,Doukhan2013corr} proved that the model (\ref{modYtilde}) has a
 strictly stationary solution $(\tilde{Y}_t, \tilde{\lambda}_t )_{t \in \Z}$ which is $\tau$-weakly dependent with finite moment of any order
 (see also Doukhan and Wintenberger (2008)\nocite{Doukhan2008}).
 But, such result cannot be apply to process $Y$ satisfying (\ref{NonStmodel_par}), since the structure of the past before $k^*$ is unknown.
 The following proposition shows that if $(Y_t)_{t \leq k^*}$ has finite moment of any order, then it also holds for $(Y_t)_{t > k^*}$.
 \begin{Prop} \label{prop1}
 Assume \textbf{A}$_0 (\Theta)$ holds. Let $Y=(Y_t)_{t \in \Z}$ satisfying (\ref{NonStmodel_par}). For any $r\geq 1$, if there exists $C_{r,0}$ such that
 $\E Y_t^r \leq C_{r,0}$ for all $t \leq k^*$, then  there exists $C > 0$ such that
 \[ \E Y_{k^* + \ell}^r \leq C ~ ~ \text{for all } \ell \geq 1 .\]
 \end{Prop}
 As we state above, $(Y_t)_{t \leq k^*}$ could be any integer-valued time series and we assume in the sequel that :
  \begin{equation}\label{moment_past}
  \text{for any } r \geq 1, \text{ there exists } C_{r,0} >0 \text{ such that } \E Y_t^r \leq C_{r,0} \text{ for all } t \leq k^* .
 \end{equation}
 \medskip
 \noindent The conditions A$_1(\Theta)$, A$_2(\Theta)$ as well as the following assumptions D$(\Theta)$, Id($\Theta$) and  Var($\Theta$) are classical for inference on such model
 see \cite{Doukhan2015}.

 \medskip
  \noindent {\bf Assumption D$(\Theta)$:} $\exists\underline{c}>0$ such that
$\displaystyle \inf_{\theta \in
 \Theta}(f_\theta(y))\geq \underline{c}$ for all $y\in (\R^+)^{\N}.$ \\
~\\
 {\bf Assumption Id($\Theta$):} For all  $(\theta,\theta')\in \Theta^2$,
 $ \Big( f_{\theta}(Y_{t-1},\dots)=f_{\theta'}(Y_{t-1},\dots)  \ \text{a.s.} ~ \text{ for some } t > k^* \Big) \Rightarrow \ \theta = \theta'.$\\
 {\bf Assumption Var($\Theta$):} For all $\theta  \in \Theta $ and $t > k^*$,
 the components of the vector $\dfrac{\partial f_{\theta}}{\partial \theta}(Y_{t-1,\dots})$  are  a.s. linearly independent. \\
  \noindent Also, we will assume in the sequel that the true parameter $\theta^*$ belongs to $\overset{\circ}{\Theta}$ (the interior of $\Theta$).
  %

  \subsection{Examples}

  \subsubsection{Linear Poisson autoregression}

   We consider an integer-valued time series $(Y_t)_{t \in \Z}$ satisfying for any $t \in \Z$
   \begin{equation}\label{Ex_linear}
   Y_t/\mathcal{F}_{t-1}  \sim \mbox{ Poisson}(\lambda_t) ~ \text{with} ~ \lambda_t = \phi_0(\theta^*) + \sum_{k \geq 1} \phi_k(\theta^*) Y_{t-k}
   \end{equation}
 with $\theta^* \in \Theta \subset \R ^d$, where the functions $\theta \mapsto \phi_k(\theta)$ are positive, twice continuous differentiable such that
  $\sum_{k \geq 1} \|\phi_k(\theta)\|_\Theta < 1$,
  $\sum_{k \geq 1} \|\phi_k '(\theta) \|_\Theta < \infty$, $\sum_{k \geq 1} \|\phi_k '' (\theta)\|_\Theta < \infty$ and $ \underset{\theta \in \Theta} {\inf} \phi_0(\theta) >0$
 (see also \cite{Doukhan2015}). Thus  Assumptions A$_i(\Theta)$, $i=0,1,2$ and D$(\Theta)$ hold.
 %
 Moreover, if there exists a finite subset $I \subset \N - \{0\}$ such that the function $\theta \mapsto (\phi_k(\theta))_{ k \in I}$ is injective, then assumption
 Id$(\Theta)$ holds and the model (\ref{Ex_linear}) is identifiable.
 Finally, assumption Var($\Theta$) holds if for any $\theta \in \Theta$, there exists $d$ functions
 $\phi_{k_1},\cdots,\phi_{k_d}$ such that the matrix $\Big(\dfrac{\partial \phi_{k_j} }{\partial \theta} \Big)_{1\leq j \leq d}$
 (computed at $\theta$) has a full rank.
 This is the case in the classical useful situations, such as for instance, the INGARCH($p,q$) model below.

  \medskip

 The classical Poisson INGARCH($p,q$) (see \cite{Ferland2006} or \cite{Weiss2009}) is obtained with
 \begin{equation}\label{INGARCH}
  \lambda_t = \alpha^*_0  + \sum_{k = 1}^{p} \alpha^*_k \lambda_{t-k} + \sum_{k = 1}^{q} \beta^*_k Y_{t-k} ;
 \end{equation}
 the true parameter $\theta^* = (\alpha^*_0, \alpha^*_1,\cdots,\alpha^*_p, \beta^*_1,\cdots,\beta^*_q) \in \Theta$ where $\Theta$ is a compact subset of
 $(0, +\infty) \times [0, +\infty)^{p+q}$ such that  $ \sum_{k = 1}^{p} \alpha_k + \sum_{k = 1}^{q} \beta_k <1 $ for all
 $\theta = (\alpha_0, \alpha_1,\cdots,\alpha_p, \beta_1,\cdots,\beta_q) \in \Theta$.
 This model is a special case of the model (\ref{Ex_linear}) since we can find a sequence of functions $(\psi_k(\theta))_{k \geq 0}$ such that
 $\lambda_t =   \psi_0(\theta^*_0) + \sum_{k \geq 1} \psi_k(\theta^*_0)Y_{t-k} $.
%

%
 In the model (\ref{Ex_linear}), it is often holds that
 \begin{equation}\label{Ex_linear_brk}
   \lambda_t = \begin{cases}
                \phi_0(\theta^*_0) + \sum_{k \geq 1} \phi_k(\theta^*_0) Y_{t-k}  & \text{for } t \leq  k^* \\
                 \phi_0(\theta^*_1) + \sum_{k \geq 1} \phi_k(\theta^*_1) Y_{t-k}   &  \text{for } t  > k^*
  \end{cases}
 \end{equation}
 with $\theta^*_0 \neq \theta^*_1$. There exists several references in the literature (see for instance Doukhan and Kengne \cite{Doukhan2015}, Ahmad and Francq (2016)) that address the
 inference on $\theta^*_0$ based on the observations of the stationary process $(Y_t)_{t\leq k^*}$.
 These results which are heavily based on the stationarity of the process cannot  work for $\theta^*_1$.
 Section \ref{like_inf} focusses on the estimation of $\theta^*_1$ based on the nonstationary process $(Y_t)_{t > k^*}$.

     \subsubsection{Threshold Poisson autoregression}
   We consider a threshold Poisson autoregressive model defined by :
   \begin{equation}\label{Ex_threshold}
   Y_t/\mathcal{F}_{t-1}  \sim \mbox{ Poisson}(\lambda_t) ~ \text{with}
   ~ \lambda_t = \phi_0(\theta^*) + \sum_{k \geq 1} \Big( \phi_k^+(\theta^*)\max(Y_{t-k}-\ell,0) + \phi_k^-(\theta^*)\min(Y_{t-k},\ell)\Big)
   \end{equation}
 where $\phi_0(\theta)>0$, $\phi^+_k(\theta), \phi^-_k(\theta) \geq 0$ for all $\theta \in \Theta$ and $\ell \in \N$. We can also write
  \[ \lambda_t = \phi_0(\theta^*) + \sum_{k \geq 1} \Big(\phi_k^-(\theta^*)Y_{t-k} + \big(\phi_k^+(\theta^*) - \phi_k^-(\theta^*)\big)\max(Y_{t-k}-\ell,0)\Big)  .\]
 This is an example of nonlinear model called an integer-valued threshold ARCH (or INTARCH) see \cite{Doukhan2015}; see also  \cite{Franke2012} for INTARCH($1$) model.
 Such model is often used to capture piecewise phenomenon. $\ell$ is the threshold parameter of the model.
  %
 %
  If the functions $\theta \mapsto \phi_k^+(\theta) $ and $\theta \mapsto \phi_k^-(\theta) $ are twice continuously differentiable such that
  $ \sum_{k \geq 1} \max\big( \| \phi_k^+(\theta) \|_\Theta, \| \phi_k^-(\theta) \|_\Theta \big) < 1 $,
  $ \sum_{k \geq 1} \max\big( \| \frac{\partial}{\partial \theta}\phi_k^+(\theta) \|_\Theta, \| \frac{\partial}{\partial \theta}\phi_k^-(\theta)\|_\Theta ,
  \| \frac{\partial^2}{\partial \theta^2}\phi_k^+(\theta) \|_\Theta,
  \| \frac{\partial^2}{\partial \theta^2}\phi_k^-(\theta)\|_\Theta \big) < \infty $,
   then  A$_i(\Theta)$ $i=0,1,2$ hold.
  Furthermore,  Conditions on D$(\Theta)$, Id$(\Theta)$ and Var($\Theta$) are obtained as above.

 \section{Likelihood inference}\label{like_inf}
 We focus on the inference for the model (\ref{NonStmodel_par}); that is, we consider the process  $Y=(Y_t)_{t \in \Z}$ satisfying
  \begin{equation}\label{NonStmodel_par_inf}
      Y_t | \mathcal{F}_{t-1}   \sim\mbox{ Poisson}(\lambda_t)  ~ \text{with} ~   \lambda_t = \E(Y_t | \mathcal{F}_{t-1}) = f_{\theta^*}(Y_{t-1}, Y_{t-2},\cdots)    \text{ for all } t >  k^*.
   \end{equation}
 Assume that a trajectory $(Y_{k^*+1}, \dots, Y_{k^*+n})$ of the process $(Y_t)_{t > k^*}$ is observed.
  Without loss of generality, for simplifying notation, we set $k^* =0$ in this section.
 The conditional (log)-likelihood (up to a constant)  computed on a segment $T \subset \{k^*+1, k^*+2, \cdots \} $ is given by
 \[ L_n(T,\theta) = \sum_{t \in T}(Y_t\log \lambda_t(\theta)- \lambda_t(\theta)) = \sum_{t\in T} \ell_t(\theta) \text{ with }
  \ell_t(\theta) = Y_t\log \lambda_t(\theta)- \lambda_t(\theta)\]
  where $ \lambda_t(\theta) = f_{\theta}(Y_{t-1},\dots)$. In the sequel, we use the notation $f^t_{\theta} := f_{\theta}(Y_{t-1}, \ldots)$.
 An approximation of the conditional (log)-likelihood is
 \begin{equation}\label{loglik}
 \widehat{L}_n(T,\theta) = \sum_{t \in T}(Y_t\log\widehat{\lambda}_t(\theta)-\widehat{\lambda}_t(\theta)) = \sum_{t\in T}\widehat{\ell}_t(\theta) \esp{with}
 \widehat{\ell}_t(\theta) = Y_t\log\widehat{\lambda}_t(\theta)-\widehat{\lambda}_t(\theta)
 \end{equation}
 where
 $ \widehat{\lambda}_t (\theta) := \widehat{f}^t_{\theta}:= f_{\theta}(Y_{t-1}, \dots, Y_1,0,\dots)$.
 The MLE of $\theta^*$ computed on $T$ is defined by
 \begin{equation}\label{emv}
  \widehat{\theta}(T) = \argmax_{\theta \in \Theta} (\widehat{L}_n(T,\theta)).
  \end{equation}
 For any $k,k' \in \Z$ such as $k\leq k'$, denote
 \[T_{k,k'} = \{k,k+1, \ldots,k' \} .\]
 The following theorem establishes that the MLE of $\theta^*$ based on the nonstationary process $Y$ is consistent.

 \begin{Theo}\label{theo_consistent}
 Assume  $\theta^* \in  \ring{\Theta}$, $D(\Theta)$, $\mathrm{Id}(\Theta)$, $A_0(\Theta)$ and (\ref{moment_past}) hold with
 \begin{equation} \label{theo_consistent_eq}
   \alpha_j^{(0)} = O(j^{-\gamma}), ~  \text{for some }  ~ \gamma > 3/2 .
 \end{equation}
 Then, it holds that
 \[\widehat{\theta}(T_{1,n}) \xrightarrow [n\to +\infty]{a.s.} \theta^*.\]
 \end{Theo}
 To address the asymptotic normality, set
 \begin{equation}\label{Sigma_AN}
 \widetilde{\Sigma}=E\Big(\frac{1}{\tilde{f}_{\theta^*}^0}(\frac{\partial}{\partial \theta} \tilde{f}^0_{\theta^*})(\frac{\partial}{\partial \theta} \tilde{f}^0_{\theta^*})'\Big)
 \end{equation}
 where  $\tilde{f}_{\theta}^t$ is defined in (\ref{ft_Ytilte}) and $'$ denotes the transpose.
 This matrix is symmetric and positive definite  (see \cite{Doukhan2015}).
 According to the proof of Theorem \ref{theo_AN}, the matrix
  \[ \widehat{\Sigma}_n = \Big(\frac{1}{n} \sum_{t=1}^{n}\frac{1}{\widehat{f}_{\theta}^t}\big(\frac{\partial}{\partial\theta} \widehat{f}_{\theta}^t \big)
                \big(\frac{\partial}{\partial\theta} \widehat{f}_{\theta}^t \big)'\Big)\Big|_{\theta = \widehat{\theta}(T_{1,n})}  \]
 is a consistent estimator of $\Sigma$.
 The asymptotic normality of the MLE is displayed in the following theorem.
  \begin{Theo}\label{theo_AN}
  Under the assumptions of Theorem \ref{theo_consistent} and Var($\Theta$) if $A_i(\Theta)$ $i=1,2$ hold with
 \begin{equation} \label{theo_AN_eq}
  \alpha_j^{(i)} = O(j^{-\gamma}), ~  \text{for some }  ~ \gamma > 3/2 ,
  \end{equation}
 then
 \[ \sqrt{n}(\widehat{\theta}(T_{1,n}) - \theta^*) \xrightarrow [n\to +\infty]{\mathcal{D}} \mathcal{N}(0,\widetilde{\Sigma}^{-1} )\]
 \end{Theo}

  \section{Sequential change-point detection}
  Let $(X_1,\cdots,X_m)$ be the historical observations generated according to (\ref{model2}) with the parameter $\theta^*_0$.
  We focus on the online change-point detection in the model (\ref{model}) and consider the following hypothesis testing :
  \medskip

 \noindent  $\mathbf{H_0}$: $\theta^*_0$ is constant over the observations $Y_1,\cdots,Y_m,Y_{m+1},\cdots$ {\it i.e.}  $(Y_t)_{t\in \N}$
             satisfying (\ref{model}) with  $\theta^*_0 = \theta^*_1$ ;

 \medskip
 \noindent $\mathbf{H_1}$ : the process $(Y_t)_{t\in \N}$ satisfying (\ref{model}) with  $\theta^*_0 \neq \theta^*_1$.

 \medskip

\noindent The MLE of $\theta^*_0$, computed on the historical observations is defined by
 \begin{equation}\label{emv_hist}
  \widehat{\theta}(T_{1,m})= \underset{\theta \in \Theta}{\argmax} (\widehat{L}(T_{1,m},\theta)).
  \end{equation}
  According to Section \ref{like_inf} (see also \cite{Doukhan2015}), this estimator is consistent and asymptotically normal.
 The asymptotic covariance matrix of $\widehat{\theta}(T_{1,m})$ is $\Sigma^{-1}$ with
 \begin{equation}\label{Sigma}
 \Sigma=E\Big(\frac{1}{f_{\theta_0^*}^0}(\frac{\partial}{\partial \theta} f^0_{\theta_0^*})(\frac{\partial}{\partial \theta} f^0_{\theta_0^*})'\Big).
 \end{equation}

  \medskip

  Recall that, the fluctuation-type test proposed by Chu {\it et al.} \cite{Chu1996} is based on the discrepancy between the estimators of the model's parameters.
  The classical idea of the fluctuation test is to evaluate at the monitoring step $m+k$, the distance between $\widehat{\theta}( T_{1,m} )$ and
  $\widehat{\theta}( T_{1, m+k} )$ ; by expecting this will be large enough if a change occurs at time $m+k^*$ (with $k^* < k$).
  Such idea has been employed by Na {\it et al.} \cite{Na2011}, Kengne \cite{Kengne2015}, among others.
  As pointed out by Bardet and Kengne \cite{Bardet2014}, the recursive estimator $\widehat{\theta}( T_{1,m+k})$  heavily depends on the historical data and the detection
  delay of such procedure may not be quite efficient.

 \medskip

  We follow the ideas of Bardet and Kengne \cite{Bardet2014} and propose a procedure which is based on the detector :

  \[ \widehat{C}_{k,\ell}:=  \sqrt{n}\, \dfrac{k-\ell}{k}\, \big\| \widehat{\Sigma}_m^{-1/2}   \big(\widehat{\theta}(T_{\ell,k}) - \widehat{\theta}(T_{1,m})\big)\big\| \]
 defined for any $k>m$ and $\ell=n,\cdots,k$ ; where
 \[ \widehat{\Sigma}_m = \Big(\frac{1}{m} \sum_{t=1}^{m}\frac{1}{\widehat{f}_{\theta}^t}\big(\frac{\partial}{\partial\theta} \widehat{f}_{\theta}^t \big)
                \big(\frac{\partial}{\partial\theta} \widehat{f}_{\theta}^t \big)'\Big)\Big|_{\theta =  \widehat{\theta}(T_{1,m}) }  ; \]
 is a consistent estimator of $\Sigma$ (see Section \ref{like_inf} and also \cite{Doukhan2015}). $\widehat{\Sigma}_m$ is also asymptotically symmetric and positive definite, and the detector
 $\widehat{C}_{k,\ell}$ is well defined for $m$ large enough.

   \medskip

   To avoid some distortion in the computation of $\widehat{\theta}(T_{\ell,k})$ (when $\ell$ is close to $k$),
   we introduce a  sequence of integer numbers $(v_m)_{m \in \N}$ with $v_m<< m$ and compute $\widehat{C}_{k,\ell}$ for $\ell \in \{ m-v_m, m-v_m+1,\cdots,k-v_m\}$.
   Thus, for any $k>m$ denote
    \[  \Pi_{m,k} :=   \{ m-v_m, m-v_m+1,\cdots,k-v_m\} . \]
    For technical consideration, assume that,
    \[ v_m \to \infty \quad \mbox{and}\quad v_m/{\sqrt{m}} \to 0 ~~(m\to \infty). \]
    %


   \medskip
   \noindent Note that, for any $\ell \in  \Pi_{m,k}$ both $\widehat{\theta}(T_{\ell,k})$ and $\widehat{\theta}(T_{1,m})$ are estimator of $\theta^*_0$ if change
   does not occur at time $k>m$, they are asymptotically  close and the detector $\widehat{C}_{k,\ell}$ is not too large under $H_0$.\\
     Let $T>1$ ($T$ can be equal to infinity). The monitoring scheme rejects $H_0$ at the first time $k$ satisfying  $m<k\leq [Tm]+1$ and there exists $\ell \in \Pi_{m,k}$
     such that  $\widehat{C}_{k,\ell}>c$ for a  suitably chosen constant $c>0$, where  $[x]$ denotes the integer part of $x$.

 \medskip
 \noindent To be more general, we will use a function $b:(0,\infty)\mapsto (0,\infty)$, called a boundary function satisfying:

  \medskip
 \noindent {\bf Assumption B:}  $b:(0,\infty)\mapsto (0,\infty)$ is a
  non-increasing and continuous function such that $ \underset{0<t<\infty} {\mbox{Inf}}b(t)>0$.

  \medskip
  \noindent Then the monitoring scheme rejects $H_0$ at the first time $k$ (with $n< k \leq [Tm] +1$) such that there exists $\ell \in \Pi_{m,k}$ satisfying  $\widehat{C}_{k,\ell}>b((k-\ell)/n)$.
  Hence, define the stopping time:
  \begin{align*} 
   \tau(m)&:= \text{Inf}\Big  \{ { m< k < [Tm] +1} ~ \big{/} ~ \exists \ell \in \Pi_{m,k}, ~ \widehat{C}_{k,\ell}>b((k-\ell)/m)\Big  \} \\
   &= \text{Inf} \Big\{ {m< k < [Tm] +1} ~ \big{/}  ~   \underset{\ell \in \Pi_{m,k}} {\mbox{max}} \dfrac{\widehat{C}_{k,\ell}}{b((k-\ell)/m)}>1 \Big\}  
 \end{align*} 
 with the convention that $\text{Inf}\{\emptyset\} = \infty$.  Therefore, we have
 \begin{align}
\nonumber  P\{ \tau(m)<\infty \} &=P\Big\{ ~ \underset{\ell \in \Pi_{m,k}} {\mbox{max}} \dfrac{\widehat{C}_{k,\ell}}{b((k-\ell)/m)}>1 ~ \text{ for some }
{ k ~ \text{ between } m \text{ and } [Tm] +1} \Big \} \\
 \label{prob_stop_tim} & = P\Big\{ ~ \underset{{m< k < [Tm] +1} } {\mbox{sup}} ~ \underset{\ell \in \Pi_{m,k}} {\mbox{max}} \dfrac{\widehat{C}_{k,\ell}}{b((k-\ell)/m)}>1 \Big \} .
  \end{align}
  The challenge is to choose a suitable boundary function $b(\cdot)$ such that for some given $\alpha \in (0,1)$,
  \[ \lim_{m \rightarrow \infty} P_{H_0}\{ \tau(m)<\infty \}=\alpha \]
  and
  \[ \lim_{m \rightarrow \infty}  P_{H_1}\{ \tau(m)<\infty \}=1 \]
 where the hypothesis $H_0$ and $H_1$ are  formulated above. \\
In the case where $b(\cdot)$ is a constant positive value,
 $b \equiv c$ with $c>0$, these conditions lead to compute a threshold $c=c_\alpha$ depending on $\alpha$.
 If change is detected under $H_1$ {\it i.e.}  $ \tau(m)<\infty $ and $\tau(m)>k^*$, then the detection delay is defined by
 \begin{equation}\label{dm}
 \widehat{d}_m= \tau(m) - k^*.
 \end{equation}
 $\widehat{d}_m$ is used to assess the efficiency of the procedure to early detect changes in the model. The smaller is the detection delay, the better is the efficiency
  under the alternative.

 \medskip
%

 \subsection{Asymptotic under the null hypothesis}
  Under $H_0$, all the observations are generated from the model (\ref{model2}) according to the parameter $\theta^*_0$.
  The following theorem displays the asymptotic behavior under the null hypothesis of the detector $\widehat{C}_{m,k}$ for the open and closed-end procedure.
  \begin{Theo}\label{theo1}
 Assume  $D(\Theta)$, $\mathrm{Id}(\Theta)$,  Var($\Theta$) and $A_i(\Theta)$ $i=0,1,2$ hold with
 \[ \alpha_j^{(i)} = O(j^{-\gamma}), ~  \text{for some }  ~ \gamma > 3/2 .\]
 Under H$_0$ with $\theta_0^* \in  \ring{\Theta}$, for the open-end ($T=\infty$) and closed-end ($T<\infty$) procedure it holds that
      \begin{equation} \label{theo1_ep}
     \lim_{m \rightarrow \infty} P\{ \tau(m)  <\infty \} = P\Big\{\sup_ {1<t\leq T} \sup_{1<s<t} \dfrac{\| W_d(s) - s W_d (1))\|}{t~b(s)}>1 \Big\}. ,
      \end{equation}
 where $W_d$ is a $d$-dimensional standard Brownian motion.
 \end{Theo}
 Assume that $b(s) = c b_0(s)$ for some function $b_0$ satisfying the assumption \textbf{B}, with $c>0$. Thus, at a nominal level $\alpha \in (0,1)$,
 the monitoring procedure stops and rejects $H_0$ at the first time $k$ (with $1< k \leq [Tm] +1$  ) such that
 \[ \underset{\ell \in \Pi_{n,k}}{\mbox{max}} \dfrac{\widehat{C}_{k,\ell}}{b_0((k-\ell)/n)}   > c_{\alpha} \]
 where
 $c_{\alpha}$ is the $(1-\alpha)$-quantile of the distribution of $\sup_ {1<t\leq T} \sup_{1<s<t} \dfrac{\| W_d(s) - s W_d (1))\|}{t~b_0(s)} $.

 \medskip

 In Section \ref{num_result}, we will use the most "natural" boundary function $b(\cdot)=c$ where $c$ is a positive constant.
 In this case, it follows directly from Theorem \ref{theo1} that
 \[   \lim_{n \rightarrow \infty} P\{\tau(n)<\infty \} = P\{  U_{d,T} > c \} \]
 where
 \begin{equation} \label{U_d_T}
   U_{d,T} = \sup_ {1<t\leq T} \sup_{1<s<t} \frac{1}{t} \| W_d(s) - s W_d (1))\|.
 \end{equation}
 Proposition 4.1 of Bardet and Kengne \cite{Bardet2014} provides a way to compute the quantile of the distribution of $U_{d,T}$, from which the critical value of the test can be obtained.

 \subsection{Asymptotic under the alternative}
 Under the alternative, a change occurs at time $k^*> m$ and contrary to some recent works (for instance: Franke {\it et al.} \cite{Franke2012},   Doukhan and Kengne \cite{Doukhan2015}, Kengne \cite{Kengne2015}, Kirch and Tadjuidje Kamgaing \cite{Kirch2015}, Diop and Kengne \cite{Kengne2017}, Kirch and Weber \cite{Kirch2018}, $\cdots$), we do not set any additional assumption on the process after the change-point.
 Many recent works impose stationarity after the change-point. This assumption is too strong for autoregressive process ; note that, in model (\ref{model}) with $t> k^*$,
 \[ \lambda_t =  f_{\theta_1^*}(Y_{t-1}, Y_{t-2},\cdots)  \]
 depends on $\theta_1^*$ and it is contaminated by observations which depends on $\theta_0^*$.
 This shows that,  stationarity assumption on the observations after change-point is quite questionable; and that, nonstationary approach seems to be suitable.
 The proof of the following theorem is heavily based on the result of Theorem \ref{theo_consistent}.
 The results below show that the proposed monitoring procedure is consistent under the alternative for both the open-end and the closed-end methods.
 %
%
 %
 \begin{Theo}\label{theo2}
 Assume  $D(\Theta)$, $\mathrm{Id}(\Theta)$,  Var($\Theta$) and $A_i(\Theta)$ $i=0,1,2$ hold with
 \[ \alpha_j^{(i)} = O(j^{-\gamma}), ~  \text{for some }  ~ \gamma > 3/2 .\]
 Under the alternative H$_1$, if $\theta_0^*, \theta_1^* \in  \ring{\Theta}$ and there exists $T^* \in (1,T)$ such that $k^* = [T^* m]$,
  for the open-end ($T=\infty$) and closed-end ($T<\infty$) procedure, then for $k_m=k^*(m)+m^\delta $ with $ \delta \in (1/2,1)$, it holds that
      \begin{equation} \label{theo2_ep}
   \max_{\ell \in \Pi_{m,k_m}} \, \dfrac{\widehat{C}_{k_m,\ell}}{b((k_m-\ell)/m)} \limitepsn  \infty .
        \end{equation}
 \end{Theo}

 \noindent The  Corollary \ref{cor2} follows immediately from Theorem \ref{theo2}.
 \begin{cor}\label{cor2}
Under the assumptions of Theorem \ref{theo2},
  \[  \lim_{n \rightarrow \infty} P\{ \tau(m) <\infty \} = 1 . \]
 \end{cor}
 Hence,  it follows from Theorem \ref{theo2} that with probability one, the change is asymptotically detected both for open-end and closed-end (when $T^*<T$)
   procedures and the detection delay $\widehat{d}_n$ can be bounded by
    $\mathcal{O}_P(m^{1/2+\varepsilon})$ for any $\varepsilon>0$ (or even by $\mathcal{O}_P\big (\sqrt m (\log m)^a \big )$ with $a>0$ using the same kind of proof).

 \section{  Some numerical results}\label{num_result}
  In this section, we conduct a small simulation study and a real data example to display some empirical performances of the proposed sequential change-point procedure. 
 We focus on the  closed-end procedure with $T=1.5$; that is, the historical available data are $X_1,\cdots, X_m$ and the monitoring period is $\{m+1,\cdots, 1.5m\}$.
 In the sequel, the detector of the sequential procedure is computed with $v_m = m^{\delta}$ for $2\leq \delta \leq 5/2$.
 The corresponding quantile of the distribution of $U_{d,T}$ can be found in the Table 1 of Bardet and Kengne \cite{Bardet2014}. 
 \subsection{Sequential change-point detection in Poisson INGARCH }
\begin{figure}[h!]
\begin{center}
\includegraphics[height=8.01cm, width=16.01cm]{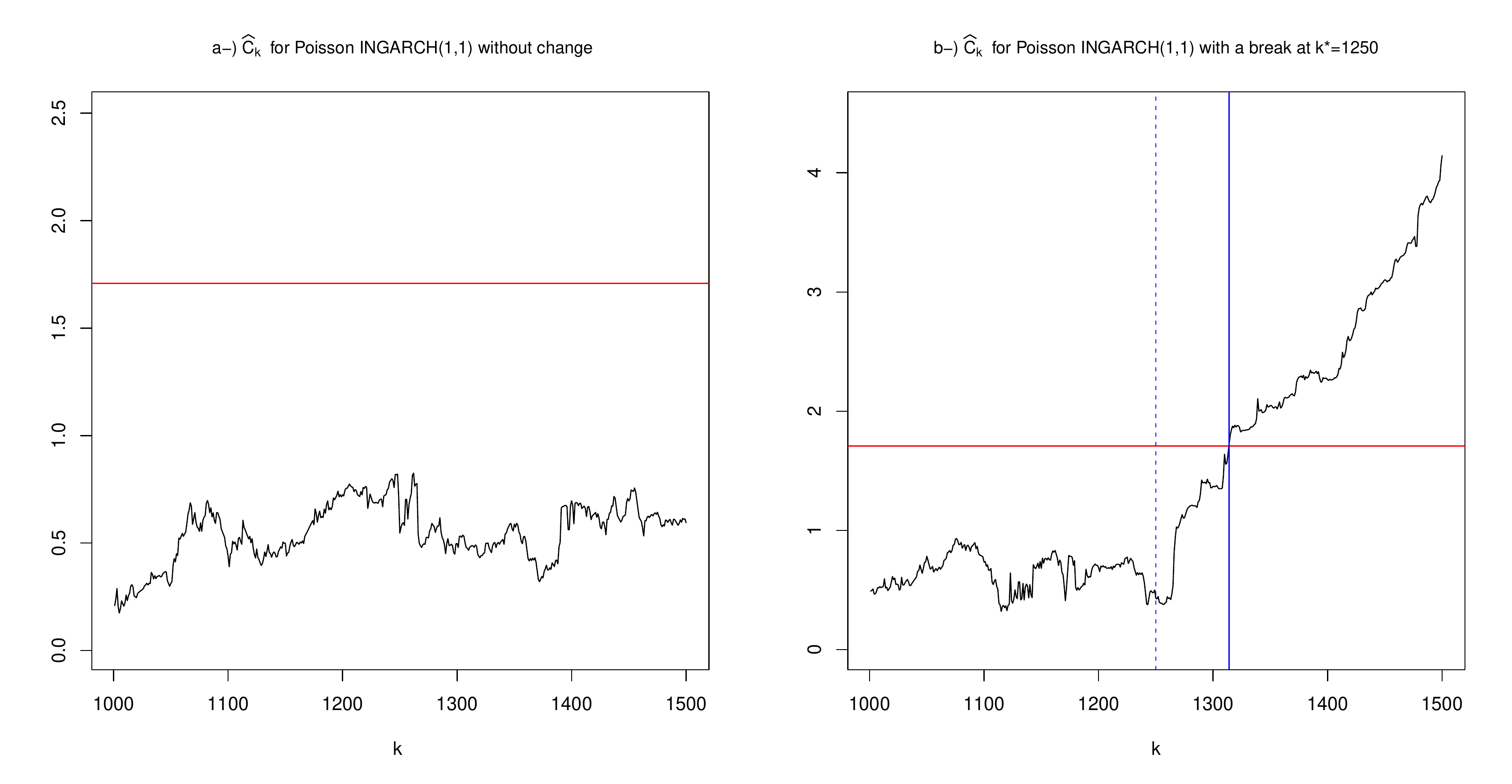}
\end{center}
\caption{ A realization of the detector  $(\widehat{C}_{k})_{1001\leq k \leq 1500}$  for a Poisson INGARCH$(1,1)$ with $m=1000$.
   a-) The parameter $\theta^*_0=(1,0.2,0.15)$ is constant;
    b-) the parameter $\theta^*_0=(1,0.2,0.15)$  changes to $\theta^*_1 = (1,0.2,0.5)$ at $k^*=1250$.
   The horizontal solid line represents the limit of the critical region, the vertical dotted line indicates where the change occurs and the vertical solid
   line indicates the time where the sequential procedure detects a  break in the observations.} \label{Online_Poisson_INGARCH}
\end{figure}

 We consider a Poisson INGARCH(1,1)
 \begin{equation}\label{simul_ingarch11}
  Y_t/\mathcal{F}_{t-1}  \sim \mbox{ Poisson}(\lambda_t) ~ \text{with} ~ \lambda_t = \alpha^*_0 + \alpha^*_1\lambda_{t-1} +  \beta^*_1Y_{t-1}
 \end{equation}
 where  $\theta^*_0 = (\alpha^*_0,\alpha^*_1,\beta^*_1)$ denote the parameter of the model.
 For any $k>m$, denote 
 $ \widehat{C}_{k} = \underset{\ell \in \Pi_{n,k}} {\mbox{max}} \widehat{C}_{k,\ell} $.
 For n=1000, Figure \ref{Online_Poisson_INGARCH} displays the statistics $(\widehat{C}_{k})_{1001\leq k \leq 1500}$ in a scenario without change a-) and a scenario with a change-point at $k^* = 1.25m=1250$ b-).
 Figure \ref{Online_Poisson_INGARCH} a-) shows that the detector $\widehat{C}_{k}$ is under the horizontal line that defined the critical region of the test; whereas in Figure 1 b-), the detector is under the horizontal before change occurs, and increases with a high rate until exceed the critical value after the change-point. As pointed out by Bardet and Kengne \cite{Bardet2014}, such growth rate over a long period indicates that something is happening in the model. 
 We consider the model (\ref{simul_ingarch11}) with scenarios under H$_0$ and H$_1$ with break at $k^* = 1.25 m$.
 Table \ref{tab1} indicates the empirical levels and powers based on 100 replications for $m=200, 500, 1000$.
    \begin{table}[htbp]
    \centering
    \begin{tabular}{c c c c c}
      \hline
           &     & $m=200$&   $m=500$ &  $m=1000$  \\
      \hline
      \hline
   Empirical levels :  &  $\theta^*_0=(1,0.2,0.15)$ & 0.08 & 0.06 & 0.05 \\
                        &  $\theta^*_0=(0.75,0.5,0.3)$ & 0.09 & 0.07 & 0.06 \\ 
                       &  $\theta^*_0=(2.5,0,0.35)$ & 0.07 & 0.06 & 0.04 \\
    &   &       &     &        \\
     Empirical powers :  &  $\theta^*_0=(1,0.2,0.15)$ ; ~ $\theta^*_1 = (1,0.2,0.5)$  & 0.63  & 0.96   &  0.98 \\
     & $\theta^*_0=(0.75,0.5,0.3)$; ~ $\theta^*_1=(0.25,0.5,0.3) $ & 0.60 & 0.92 & 0.97 \\
                          & $\theta^*_0=(2.5,0,0.35)$; ~ $\theta^*_1=(4.5,0.05,0.6) $ & 0.81  & 1  &   1 \\
      \hline
    \end{tabular}
 \caption{ Empirical levels and powers  for sequential change-point detection in Poisson INGARCH($1,1$) model. The empirical levels are computed when  $\theta_0^*=(1,0.2,0.15), (0.75,0.5,0.3), (2.5,0,0.35)$ is constant (under $H_0$) and the empirical powers when $\theta_0^*=(1,0.2,0.15), (0.75,0.5,0.3), (2.5,0,0,0.35)$ changes respectively to $\theta_1^*=(1,0.2,0.5), (0.25,0.5,0.3), (4.5,0.05,0.6)$ (under the alternative) at $k^*=1.25m$. }\label{tab1}
 \end{table}
 Some elementary statistics of the empirical detection delay (defined at (\ref{dm})) are summarized in Table \ref{tab2}.
   \begin{table}[htbp]
    \centering
    \begin{tabular}{ c c c c c c c c c}
      \hline
          $\widehat{d}_n$   &   &     Mean  & SD  &   Min  &  $Q_1$  &  Med  &  $Q_3$ & Max \\
                     \hline
                      \hline 
  $\theta^*_0=(1,0.2,0.15)$ ; ~ $\theta^*_1 = (1,0.2,0.5)$  & $m=200$ ;  $k^*=250$ ~ & 34.92 & 10.42 & 11 & 27 & 36 & 43  & 50  \\ 
         & $m=500$ ; $k^*=625$ ~ & 59.96 & 21.94 & 17 & 45 & 61 & 72& 119 \\
      & $m=1000$ ;  $k^*=1250$ ~ & 86.19 & 33.82 & 24 &56& 85& 110 & 168 \\  
                   &                  &            &   &   &   &       &      &    \\
  $\theta^*_0=(0.75,0.5,0.3$ ; ~ $\theta^*_1 = (0.25,0.5,0.3)$  & $m=200$ ;  $k^*=250$ ~ & 32.62    & 14.05   & 3   &  23  & 36  &  44 & 50 \\ 
         & $m=500$ ; $k^*=625$ ~ & 71.82 &  20.17  & 13   &  59 &  80  & 87 & 98 \\
      & $m=1000$ ;  $k^*=1250$ ~ & 103.5 & 36.53  & 15  & 96 & 109 &  116  &  183  \\  
                   &                  &            &   &   &   &       &      &    \\ 
 $\theta^*_0=(2.5,0,0.35)$; ~ $\theta^*_1=(4.5,0.05,0.6) $     & $m=200$ ;  $k^*=250$~ & 27.78 & 12.03 & 6 & 19 & 26 & 39 & 48 \\ 
      & $m=500$ ;  $k^*=625$  ~ & 55.38 & 24.01 & 5 & 34 & 56 & 77 & 91 \\
   & $m=1000$ ; $k^*=1250$ ~ & 56.22 & 36.43 & 7 & 25 & 42 & 97 & 120 \\  
      \hline
    \end{tabular}
 \caption{{\footnotesize  Elementary statistics of the empirical detection delay for sequential change-point detection change in a Poisson INGARCH(1,1).}}
        \label{tab2}
\end{table}

 \medskip
The results of Table \ref{tab1} displays some distortion in the empirical levels for the first scenario when $n=200$ and the second scenario when $n = 200, 500$. But the empirical level decreases as $n$ increases and for the three cases, it is close to the nominal level for $1000$. Also, empirical powers
increase with $n$ and approaching one when $n = 1000$ for the three scenarios. These results are consistent with Theorem \ref{theo1} and  Corollary \ref{cor2}. 
In Table \ref{tab2}, for example, when $n = 200$ with the break occurred at the time $k^* = 250$, this break is detected on average after a delay of 35, 33 and 28 respectively for these scenarios.
One can see that, for two historical sample sizes $m_1$ and $m_2$ with $m_1 < m_2$, the sequence $\widehat{d}_{m_2} - \sqrt{m_2/m_1} \widehat{d}_{m_1} $ decreases when $m_1$ and $m_2$ increases and it is on average, close or less than 0 when $m_1=500$ and $m_2 = 1000$. 
This is in accordance with Theorem \ref{theo2} where $\widehat{d}_m$ can be bounded by  $\mathcal{O}_P\big (m^{1/2 + \epsilon} \big)$ for any $\epsilon>0$.
  
 \subsection{Real data example} 
\begin{figure}[h!]
    \centering
    \includegraphics[width=12.05cm,height=9.05cm]{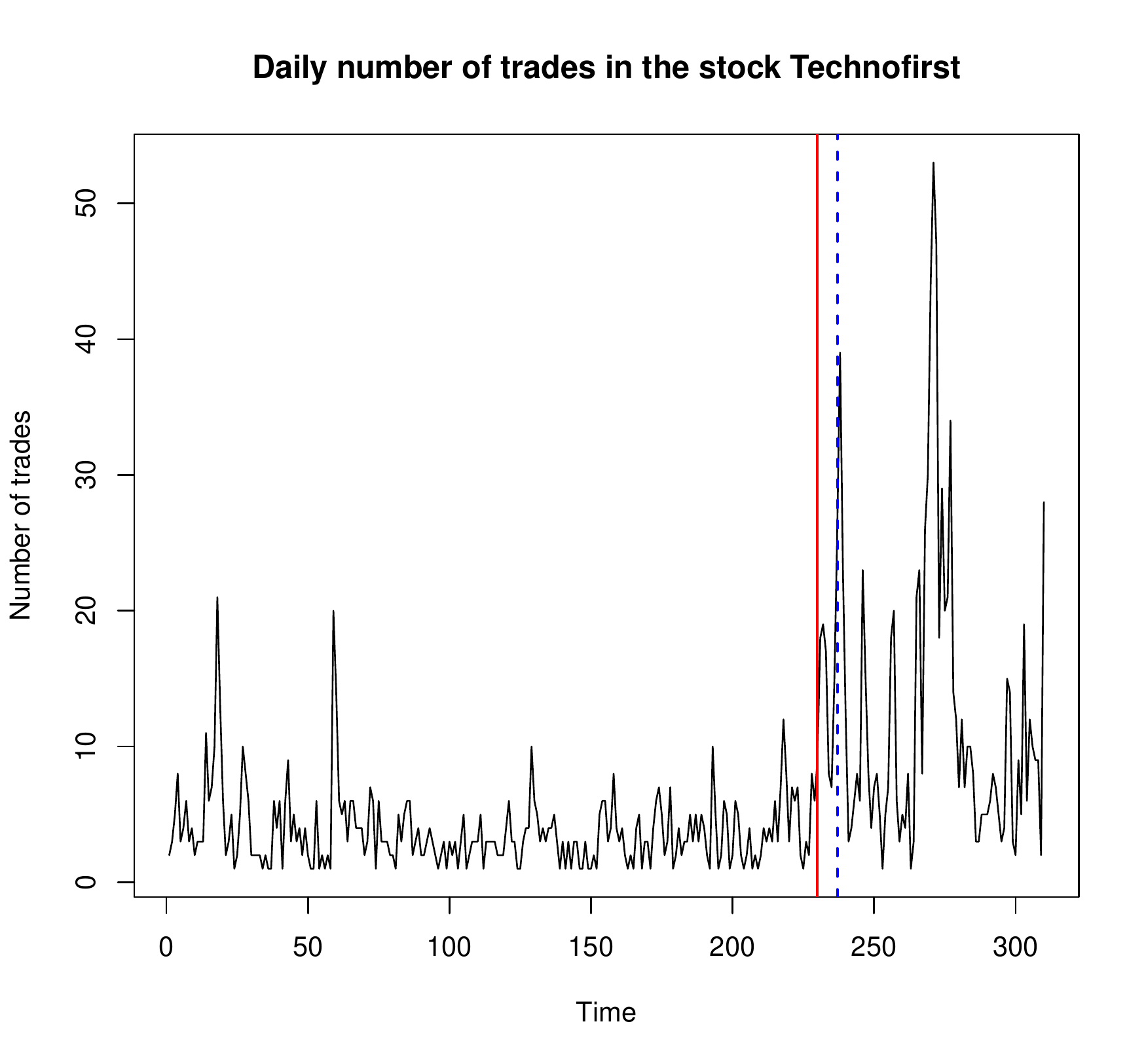}
    \caption{ Daily number of trades in the stock of Technofirst from 04 January 2010 to 05 September 2011.  The solid line represents the break that has been detected by Diop and Kengne \cite{Kengne2019} from a retrospective procedure. The dotted line   indicates the stopping time of the sequential procedure proposed.}
    \label{Technofirst}
 \end{figure}

  We consider the daily number of trades in the stock of Technofirst listed in the NYSE Euronext group.
These data have been analyzed by Ahmad and Francq \cite{Ahmad2016} with the PQMLE, and have concluded that the INGARCH(1,3) is more appropriate.
Diop and Kengne (2019) have applied the multiple change-point with an INGARCH(1,1) representation based on the Poisson quasi-maximum likelihood estimator. 
 We consider the data from 04 January 2010 to 05 September 2011 (see Figure \ref{Technofirst}); there are 310 observations. 
 For the data from $t=1$ to $t=230$, Diop and Kengne \cite{Kengne2019} have showed that the INARCH(1) representation is more appropriate and the INGARCH(1,1) representation has been used for $t>230$. So, we applied the Poisson INGARCH(1,1) model and consider the observations from $t=1$ to $t=207$ as the historical data. 
 We carry out the  sequential procedure in the closed-end setting with $T=1.5$; so, $[T \times m] = 310$. Therefore,  the monitoring starts at the time $t = 208$.
 The estimation of the parameter computed on the historical data is $\widehat{\theta}_0 = (2.43,2\times 10^{-8},0.35)$.

\medskip
 
\begin{figure}[h!]
    \centering
    \includegraphics[width=12.05cm,height=9.05cm]{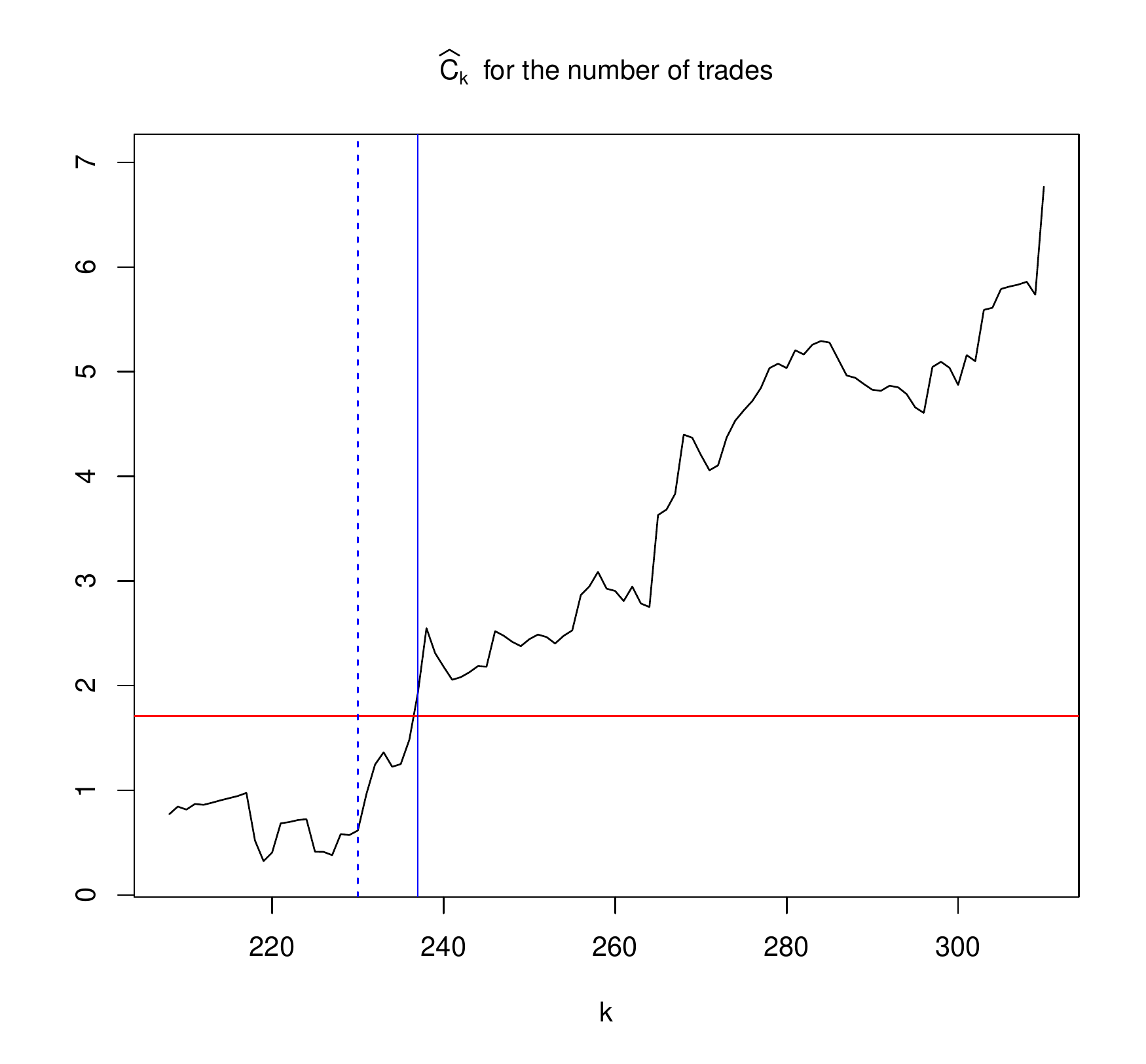}
    \caption{ Realizations of the statistics  $(\widehat{C}_{k})_{208\leq k \leq 310}$ for the daily number of trades in the stock of Technofirst from 04 January 2010 to 05 September 2011 ;
    the historical data considered are the first $207$ observations.  The horizontal solid line represents the limit of
    the critical region, the vertical dotted line represents the break that has been detected by using the retrospective procedure of Diop and Kengne \cite{Kengne2019} and the
    vertical solid line indicates the stopping time of the sequential procedure.}
    \label{Online_Technofirst}
 \end{figure}

 Figure \ref{Online_Technofirst} displays the realizations of the detector $ \widehat{C}_{k} = \underset{\ell \in \Pi_{n,k}} {\mbox{max}} \widehat{C}_{k,\ell} $, with $k=208,\cdots,310$. One can see that the sequential procedure stops at time $t=237$. In term of the detection delay, it appears that the procedure works well for this real data example; in the sense that the sequential procedure stops 7 days after the break time detected by Diop and Kengne \cite{Kengne2019}. 

 \section{Concluding remarks}
 This work addresses the question of inference for nonstationary time series of counts. After a time $k^*$, the process is a nonstationary Poisson autoregressive model with the conditional mean that depends on a parameter $\theta^*$. We carry out an approximation study between this process and the stationary regime; which allows us to establish that the MLE of $\theta^*$ computed with the nonstationary observations is consistent and asymptotically normal. We thus provide a detailed proof of an issue that has been addressed by Doukhan and Kengne \cite{Doukhan2015} (see Remark 4.1).
 These results are very useful both in retrospective and in sequential change-point problem. We perform an application to sequential change-point detection and propose a consistent procedure which the detection delay can been bounded by $\mathcal{O}_P\big (m^{1/2 + \epsilon} \big)$ for any $\epsilon>0$. Empirical studies show that the procedure works well for simulated and real data example with satisfactory detection delay.
 An extension of this work is the study of the inference for nonstationary model where the conditional distribution is different from Poisson, and could be for instance negative binomial, binary, $\cdots$

 \section{ Proofs of main results} \label{proof}
  Let $ (\psi_n)_n  $ and $ (r_m)_m  $ be sequences of random variables or vectors. Throughout this section, we use the notation
  $ \psi_m = o_P(r_m)  $ to mean :  for all $  \varepsilon  > 0, ~ P( \|\psi_m \| \geq \varepsilon \|r_m \| ) \rightarrow 0 $ as $m \rightarrow  \infty$.
  Write $ \psi_m = O_P(r_m)  $ to mean :  for all  $  \varepsilon > 0 $,  there exists $C>0$  such that  $ ~ P( \|\psi_m \| \geq C \|r_m \| )\leq \varepsilon  $
   for $n$ large enough.

 \medskip

  \paragraph{ Proof of Proposition \ref{prop1}}  ~
 We will prove that, for all $r \in \N$, there exists $C_r > 0$ such that
  \begin{equation} \label{proof_prop1_induction_r}
  \E Y^r_{t} \leq C_r, ~ \forall t \in \Z .
  \end{equation}
  \noindent Recall that for all $\ell \geq 1$,
  \[  Y_{k^* + \ell} | \mathcal{F}_{k^* + \ell-1}   \sim\mbox{ Poisson}(\lambda_{k^* + \ell}) ~ \text{with} ~
  \lambda_{k^* + \ell} =   f_{\theta^*}(Y_{t-1}, Y_{t-2},\cdots)  = f^{k^* + \ell}_{\theta^*} .\]
  According to assumption \textbf{A}$_0 (\Theta)$, we have for all $\ell \geq 1$,
 \begin{equation} \label{proof_prop1_Yk_ell_lip}
  f^{k^* + \ell}_{\theta^*}  \leq | f^{k^* + \ell}_{\theta^*}  -  f_{\theta^*}(0)  |  + f_{\theta^*}(0)
 \leq \sum_{j \geq 1} \alpha_j^{(0)}Y_{k^* + \ell -j } + f_{\theta^*}(0) .
 \end{equation}
In the sequel, we set $\alpha^{(0)} = \sum_{j \geq 1} \alpha_j^{(0)} $.
 If (\ref{proof_prop1_induction_r}) holds for some $r \in \N$, then we get from the Jensen's inequality,
  \begin{equation} \label{proof_prop1_Jen_inq}
   \E \big[ \big(  \sum_{j \geq 1} \alpha_j^{(0)}Y_{k^* + \ell -j } \big)^r  \big]
    = (\alpha^{(0)})^r \E \big[ \big(  \sum_{j \geq 1} \frac{ \alpha_j^{(0)} }{\alpha}  Y_{k^* + \ell -j } \big)^r  \big]
     \leq  (\alpha^{(0)})^{r-1}   \sum_{j \geq 1}  \alpha_j^{(0)}    \E Y_{k^* + \ell -j }^r  \leq (\alpha^{(0)})^r C_r .
   \end{equation}
 Moreover, under (\ref{proof_prop1_induction_r}), for some $r \in \N$ since $Y_t^s \leq Y_t^r$ $a.s.$ for any $s \leq r$, we have $\E Y_t^s \leq C_r$ for $s \leq r$. Thus, we can get
 $C_s \leq C_r$ for any $s \leq r$.
 Therefore, for all $\ell \geq 1$,
 \begin{align} \label{proof_prop1_fklr}
   \nonumber \E [ (f^{k^* + \ell}_{\theta^*} )^r ] &\leq \E\big[ \big( \sum_{j \geq 1} \alpha_j^{(0)}Y_{k^* + \ell -j } + f_{\theta^*}(0) \big)^r  \big] \\
  \nonumber &\leq  \sum_{ s= 0}^r  \dbinom{r}{s} \big( f_{\theta_1^*}(0) \big)^{r-s} \E\big[ \big( \sum_{j \geq 1} \alpha_j^{(0)}Y_{k^* + \ell -j } \big)^s  \big] \\
   \nonumber &\leq  \sum_{ s= 0}^r  \dbinom{r}{s} \big( f_{\theta^*}(0) \big)^{r-s} (\alpha^{(0)})^s C_s  \big] \\
   &\leq  \sum_{ s= 0}^r C_r \dbinom{r}{s} \big( f_{\theta^*}(0) \big)^{r-s} (\alpha^{(0)})^s  \leq C_r \big( \alpha^{(0)} + f_{\theta^*}(0) \big)^r \leq C_{r,f}
    \end{align}
 with $  C_{r,f} = C_r \big( \alpha^{(0)} + f_{\theta^*}(0) \big)^r  $. \\
 Let us show by induction that for all $r \in \N$, there exists $C_r > 0$ such that (\ref{proof_prop1_induction_r}) holds.
 For $r=1$, if $C_1$ exists, we will have $\E Y_t \leq C_1$ for all $t \leq k^*$; and according to (\ref{proof_prop1_Yk_ell_lip}), for all $\ell \geq 1$,
 \[ \E Y_{k^* + \ell} = \E f^{k^* + \ell}_{\theta_1^*} \leq \sum_{j \geq 1}  \alpha_j^{(0)}    \E Y_{k^* + \ell -j } + f_{\theta_1^*}(0) \leq \alpha^{(0)} C_1 + f_{\theta_1^*}(0) . \]
 Hence, (\ref{proof_prop1_induction_r}) holds with $C_1 = \max( C_{r,0}, \frac{1}{1 - \alpha^{(0)}} f_{\theta^*}(0) )$.
 Assume (\ref{proof_prop1_induction_r}) holds until $r \in \N$.
 According to Lemma 1 of \cite{Latour2006} (see also Lemma A.1. of \cite{Doukhan2012}) and (\ref{proof_prop1_fklr}), for all $\ell \geq 1$,
 \begin{align} \label{proof_prop1_Y_rp1_k_ell}
    \nonumber \E Y_{k^* + \ell}^{r+1} &= \E\big( \E( Y_{k^* + \ell}^{r+1} | \mathcal{F}_{k^*+\ell -1 } ) \big)  = \sum_{ s= 0}^{r+1}  {r+1 \brace s} \E[(f^{k^* + \ell}_{\theta^*} )^s] \\
     \nonumber &= \E[(f^{k^* + \ell}_{\theta^*} )^{r+1} + \sum_{ s= 0}^{r}  {r+1 \brace s} \E[(f^{k^* + \ell}_{\theta^*} )^s] \\
      & \leq \E[(f^{k^* + \ell}_{\theta^*} )^{r+1} + \sum_{ s= 0}^{r}  {r+1 \brace s} C_{s,f}
 \end{align}
 where for all $n, k \in \N_0$, ${n \brace k}$ denotes the Stirling numbers of the second kind that satisfies the recurrence
 ${n \brace k} = {n-1 \brace k-1} + k{n-1 \brace k} $ with ${n \brace n} = 1$ $\forall n \in \N_0$, ${n \brace 0} = 0$ $\forall n \in \N$ and ${n \brace k} = 0$ if $ k > n$.
 Hence,  if $C_{r+1}$ exists, it must satisfy $  C_{r+1}  \geq \E Y_t^{r+1}$ for all $t\leq k^*$ ; and according to (\ref{proof_prop1_Yk_ell_lip}) and (\ref{proof_prop1_Jen_inq}), we have
 \begin{align} \label{proof_prop1_f_rp1_k_ell}
    \nonumber \E [ (f^{k^* + \ell}_{\theta^*} )^{r+1} ] &\leq \E\big[ \big( \sum_{j \geq 1} \alpha_j^{(0)}Y_{k^* + \ell -j } + f_{\theta^*}(0) \big)^{r+1}  \big] \\
  \nonumber &\leq  \sum_{ s= 0}^{r+1}  \dbinom{r+1}{s} \big( f_{\theta^*}(0) \big)^{r-s+1} \E\big[ \big( \sum_{j \geq 1} \alpha_j^{(0)}Y_{k^* + \ell -j } \big)^s  \big] \\
   \nonumber &\leq \E\big[ \big( \sum_{j \geq 1} \alpha_j^{(0)}Y_{k^* + \ell -j } \big)^{r+1}  \big]
           + \sum_{ s= 0}^{r}  \dbinom{r+1}{s} \big( f_{\theta^*}(0) \big)^{r-s+1} \E\big[ \big( \sum_{j \geq 1} \alpha_j^{(0)}Y_{k^* + \ell -j } \big)^s  \big] \\
\nonumber &\leq (\alpha^{(0)})^{r+1}C_{r+1} + \sum_{ s= 0}^{r}  \dbinom{r+1}{s} \big( f_{\theta^*}(0) \big)^{r-s+1} (\alpha^{(0)})^sC_s \\
\nonumber &\leq (\alpha^{(0)})^{r+1}C_{r+1} + C_r \sum_{ s= 0}^{r}  \dbinom{r+1}{s} \big( f_{\theta^*}(0) \big)^{r-s+1} (\alpha^{(0)})^s \\
 \nonumber &\leq (\alpha^{(0)})^{r+1}C_{r+1} + C_r \big( (\alpha^{(0)} + f_{\theta^*}(0))^{r+1} - (\alpha^{(0)})^{r+1} \big).
 \end{align}
 Hence, (\ref{proof_prop1_Y_rp1_k_ell}) gives
 \[  \E Y_{k^* + \ell}^{r+1} \leq  (\alpha^{(0)})^{r+1}C_{r+1} + C_r \big( ( \alpha^{(0)} + f_{\theta^*}(0))^{r+1} - (\alpha^{(0)})^{r+1} \big) +  \sum_{ s= 0}^{r}  {r+1 \brace s} C_{s,f}  .\]
 Thus, (\ref{proof_prop1_induction_r}) holds with
  $ C_{r+1} = \max \Big( C_{r,0} ,
       \frac{ C_r \big( (\alpha^{(0)} + f_{\theta^*}(0))^{r+1} - (\alpha^{(0)})^{r+1} \big) +  \sum_{ s= 0}^{r+1}  {r+1 \brace s} C_{s,f} }{ 1 - (\alpha^{(0)})^{r+1} }   \Big)$.
  This completes the proof of the Proposition.
\Box

 ~ \\

 As stated in the Introduction, the following approximation study to the stationary regime plays a key role in the proof Theorem \ref{theo_consistent} and \ref{theo_AN}.

 \paragraph{ Approximation with stationary solutions after breakpoint}  ~

 \medskip

  Under the Lipschitz-type \textbf{A}$_0 (\Theta)$, there exists (see \cite{Doukhan2012,Doukhan2013corr}) a stationary solutions of the process after $k^*$ ; that is, there exists a
  stationary process $\tilde{Y} =(\tilde{Y}_t)_{t \in \Z}$  with finite moment of any order, satisfying :
   \begin{equation}\label{Ytilde}
   \tilde{Y}_t | \mathcal{\tilde{F}}_{t-1}   \sim\mbox{ Poisson}(\tilde{\lambda}_t) ~ \text{with} ~  \tilde{\lambda}_t=f_{\theta^*}(\tilde{Y}_{t-1}, \tilde{Y}_{t-2},\cdots) ~ ~
   \text{for } t \in \Z
     \end{equation}
 where $ \mathcal{\tilde{F}}_{t} = \sigma(\tilde{Y}_s, s\leq t)$ is the $\sigma$-field generated by the whole past of $\tilde{Y}$.\\
 \noindent For $T \subset \N $, let us consider the conditional  (log)-likelihood function (up to a constant) of this stationary regime computed on $T$ :
 %
  \begin{equation}\label{loglik_Ytilte}
   \tilde{L}(T,\theta) = \sum_{t \in T} \big( \tilde{Y}_t\log \tilde{\lambda}_t(\theta)- \tilde{\lambda}_t(\theta) \big)
   = \sum_{t \in T} \tilde{\ell}_t(\theta) ~ \text{ with } ~ \tilde{\ell}_t(\theta) =  \tilde{Y}_t\log \tilde{\lambda}_t(\theta)- \tilde{\lambda}_t(\theta)
   \end{equation}
  where $ \tilde{\lambda}_t(\theta) = f_{\theta}(\tilde{Y}_{t-1},\dots)$ ;  we will use the notation
   \begin{equation}\label{ft_Ytilte}
     \tilde{f}^t_{\theta} = f_{\theta}(\tilde{Y}_{t-1}, \ldots), ~ \text{ for all } ~ t \in \Z.
    \end{equation}
 The following lemma provides an approximation of the process $(Y_t)_{t > k^*}$ to the second stationary regime.
\begin{lem} \label{lem2}
 Consider the model (\ref{NonStmodel_par}) and assume that the conditions of Theorem \ref{theo_consistent} hold.
 There exists  $C>0$ such that for all  $ \ell \geq 1$,
\begin{equation}\label{lem2_approxYtilde}
   \E | Y_{k^* + \ell}  - \tilde{Y}_{k^* + \ell} |  \leq  C\Big( \inf_{1 \leq p \leq \ell } \big\{ (\alpha^{(0)})^{\ell/p} +  \sum_{ k \geq p} \alpha_k^{(0)} \big\} \Big)
\end{equation}
 where $ \alpha^{(0)} =  \sum_{ k \geq 1} \alpha_k^{(0)} $.
 \end{lem}

 \begin{dem}
   From the representation (\ref{modelPP}), we can write (see also Remark 4.1 of \cite{Doukhan2015}),
  \[  Y_{k^* + \ell} = N_{k^* + \ell}(\lambda_{k^* + \ell}) \text{ with } \lambda_{k^* + \ell} = f^{k^* + \ell}_{\theta^*} = f_{\theta^*}(Y_{k^* + \ell-1}, \ldots)\]
  and
  \[  \tilde{Y}_{k^* + \ell} = N_{k^* + \ell}(\tilde{\lambda}_{k^* + \ell}) \text{ with } \tilde{\lambda}_{k^* + \ell} = f^{k^* + \ell}_{\theta^*}
     = \tilde{f}_{\theta^*}(\tilde{Y}_{k^* + \ell-1}, \ldots) .\]
  Hence, we have
  \[  Y_{k^* + \ell} =F(Y_{k^* + \ell-1}, \ldots; N_{k^* + \ell}) \text{ and }  \tilde{Y}_{k^* + \ell} =F(\tilde{Y}_{k^* + \ell-1}, \ldots; N_{k^* + \ell}) \]
     where
  \[ F(y_1,y_2, \ldots; N_{k^* + \ell}) = N_{k^* + \ell}(f_{\theta^*}(y_1,y_2, \ldots))  \text{ for any } y_k \in \N, ~ k\geq 1 .\]
  Therefore,
       \begin{align} \label{approx_Y_Ytilde}
   \nonumber \E | Y_{k^* + \ell}  - \tilde{Y}_{k^* + \ell} | &= \E |F(Y_{k^* + \ell}, \ldots; N_{k^* + \ell}) - F(\tilde{Y}_{k^* + \ell}, \ldots; N_{k^* + \ell})| \\
   \nonumber  &= \E |N_{k^* + \ell}( f^{k^* + \ell}_{\theta^*} ) - N_{k^* + \ell}( \tilde{f}^{k^* + \ell}_{\theta^*} )| \\
   \nonumber  &= \E\big[ \E[ |N_{k^* + \ell}( f^{k^* + \ell}_{\theta^*} ) - N_{k^* + \ell}( \tilde{f}^{k^* + \ell}_{\theta^*} ) | ~
         | ~ \mathcal{F}_{k^* + \ell -1}, \mathcal{\tilde{F}}_{k^* + \ell -1} ] \big] \\
  \nonumber   &= \E|f^{k^* + \ell}_{\theta^*} -  \tilde{f}^{k^* + \ell}_{\theta^*} | \\
   \nonumber  &= \E|f_{\theta^*}(Y_{k^* + \ell-1}, \ldots)  -  f_{\theta^*}(\tilde{Y}_{k^* + \ell-1}, \ldots)  | \\
    & \leq \sum_{k\geq 1} \alpha^{(0)}_k \E|Y_{k^* + \ell-k}   - \tilde{Y}_{k^* + \ell-k} |
    \end{align}
  where the third equality holds since  $ |N_{k^* + \ell}( f^{k^* + \ell}_{\theta^*} ) - N_{k^* + \ell}( \tilde{f}^{k^* + \ell}_{\theta^*} ) | ~
  | ~ \mathcal{F}_{k^* + \ell -1}, \mathcal{\tilde{F}}_{k^* + \ell -1} $  can also be considered as a number of events $N_t$ that occur in the time interval
  $[0, |f^{k^* + \ell}_{\theta^*} -  \tilde{f}^{k^* + \ell}_{\theta^*} |]$. \\
  \noindent For all $\ell \in \Z$, set $ u_{\ell} := \E | Y_{k^* + \ell}  - \tilde{Y}_{k^* + \ell} |$.
         %
  According to Proposition \ref{prop1}, we can find a constant $C_1 > 0$ satisfying $\E Y_t \leq C_1$ for all $t \in \Z$. Hence, since the process $\tilde{Y}$ is stationary,
  we have for any $\ell \in \Z$,
      $ u_{\ell} \leq \E  Y_{k^* + \ell}    +  \E  \tilde{Y}_{k^* + \ell}  \leq  C_1 + \E \tilde{Y}_0 $.
      Set $C = C_1 + \E \tilde{Y}_0$.
   Let us show by induction on $\ell$ that for any $\ell \in \N$,
     \begin{equation}\label{proof_lem2_induction}
     u_{\ell} \leq  C\Big( \inf_{1 \leq p \leq \ell } \big\{ (\alpha^{(0)})^{\ell/p} +  \frac{1}{ 1 - \alpha^{(0)}} \sum_{ k \geq p} \alpha_k^{(0)} \big\} \Big)  .
     \end{equation}
    For $\ell = 1$, (\ref{proof_lem2_induction}) holds according to (\ref{approx_Y_Ytilde}).
    Assume that (\ref{proof_lem2_induction}) holds until $\ell$.
    Let $ 1 \leq p \leq \ell + 1$. From (\ref{approx_Y_Ytilde}), we have
     \begin{align}
     \nonumber  u_{\ell  + 1} &\leq   \sum_{k = 1}^{p-1} \alpha^{(0)}_k  u_{\ell - k +1}  +  \sum_{k \geq p} \alpha^{(0)}_k  u_{\ell - k  +1}  \\
      &\leq C \sum_{k = 1}^{p-1} \alpha^{(0)}_k \big( (\alpha^{(0)})^{(\ell-k+1)/p} +  \frac{1}{ 1 - \alpha^{(0)}}  \sum_{ i \geq p} \alpha_i^{(0)} \big)
    \label{proof_lem2_appl_induct}                                   +  C \sum_{k \geq p} \alpha^{(0)}_k \\
      &\leq C \sum_{k = 1}^{p-1} \alpha^{(0)}_k (\alpha^{(0)})^{(\ell -k +1)/p} +    C \frac{\alpha^{(0)}}{ 1 - \alpha^{(0)}}   \sum_{ i \geq p} \alpha_i^{(0)}
    \nonumber     +  C \sum_{k \geq p} \alpha^{(0)}_k \\
   \nonumber  &\leq C  (\alpha^{(0)})^{(\ell -(p-1) +1)/p} \alpha^{(0)}  +    C \frac{1}{ 1 - \alpha^{(0)}}   \sum_{ k \geq p} \alpha_k^{(0)}  \\
     \nonumber &\leq C \Big(   (\alpha^{(0)})^{(\ell +1)/p}   +    \frac{1}{ 1 - \alpha^{(0)}}   \sum_{ k \geq p} \alpha_k^{(0)} \Big)
     \end{align}
 Therefore (\ref{proof_lem2_induction}) holds for $\ell +1$. Thus, (\ref{lem2_approxYtilde}) holds.
 Note that, in the inequality (\ref{proof_lem2_appl_induct}), we have applied
  \[ u_{\ell-k+1} \leq C \big( (\alpha^{(0)})^{(\ell-k+1)/p} +  \frac{1}{ 1 - \alpha^{(0)}}  \sum_{ i \geq p} \alpha_i^{(0)} \big)     \]
  even when $p \geq \ell-k+1$. Indeed, we have in this case
  \[ u_{\ell-k+1} \leq C  \sum_{ i \geq 1} \alpha_i^{(0)} \leq C \alpha^{(0)} \leq  C (\alpha^{(0)})^{(\ell-k+1)/p}    .\]
  \end{dem}

   \paragraph{ Proof of Theorem \ref{theo_consistent}}  ~
 Let us prove that
   \begin{equation}\label{proof_theo_consistent_approxLtilde}
  \frac{ 1 }{n} \big \|\widehat{L}\big(T_{k^*+1,k^*+n},\theta\big)-  \tilde{L}\big(T_{k^*+1,k^*+n},\theta\big)  \big\|_{\Theta} \overset{\texttt{a.s.}}{\underset{n \to \infty}\longrightarrow} 0.
  \end{equation}
 Indeed, consider the function $\tilde{\mathcal{L}} : \theta \mapsto \E \tilde{\ell}_0(\theta)$; where $\tilde{\ell}_0$ is defined in (\ref{loglik_Ytilte}).
 From the proof of Theorem 3.1 of \cite{Doukhan2015}, we have $ \E\big(\underset{ \theta \in \Theta}{\sup |\tilde{\ell}_0 (\theta)|} \big) < \infty $,
 \[   \big \| \frac{ 1 }{n} \tilde{L}\big(T_{k^*+1,k^*+n},\theta\big) -  \tilde{\mathcal{L}}(\theta) \big\|_{\Theta} \limitepsm 0, \]
 and that the function $\tilde{\mathcal{L}}$  has a unique maximum at $\theta^*$.
 If (\ref{proof_theo_consistent_approxLtilde}) holds, we will get
 \[   \big \| \frac{ 1 }{n} \widehat{L}\big(T_{k^*+1,k^*+n},\theta\big) -  \tilde{\mathcal{L}}(\theta) \big\|_{\Theta} \limitepsm 0 ;\]
 and standard arguments can be used to conclude that $\widehat{\theta}(T_{k^*+1,k^*+n}) \limitepsn \theta^*$.
 Thus, to complete the proof of the Theorem, it suffices to prove (\ref{proof_theo_consistent_approxLtilde}).

 \medskip

 In the sequel, $C$ denotes a positive constant whom value may differ from an inequality to another.
  We have
\begin{multline*}
\frac{ 1 }{n} \big \|\widehat{L}\big(T_{k^*+1,k^*+n},\theta\big)-  \tilde{L}\big(T_{k^*+1,k^*+n},\theta\big)  \big\|_{\Theta}
\leq \frac{1}{n}\sum_{t \in T_{k^*+1,k^*+n}} \| \widehat{l}_t(\theta) - \tilde{l}_t(\theta) \|_{\Theta} \\
 \leq \frac{1}{n}\sum_{t =1}^{n} \| \widehat{l}_{k^* + t}(\theta) - \tilde{l}_{k^* + t}(\theta) \|_{\Theta} .
\end{multline*}
 Let $0<r<1$. According to Kounias and Weng (1969), it suffices to show that
 \begin{equation}\label{proof_lem2_Kounias_cond}
  \sum_{\ell \geq 1} \big( \frac{1}{\ell} \big)^r \E \big[ \| \widehat{l}_{k^* + \ell}(\theta) - \tilde{l}_{k^* + \ell}(\theta) \|^r_{\Theta}  \big]  < \infty.
 \end{equation}
 By using the inequality $|a_1 b_1-a_2 b_2| \leq |a_1| | b_1-b_2| + |b_2| | a_1-a_2| $ $\forall a_1, a_2, b_1, b_2 \in \R$, we get
 for all $ \ell \in \N$ and $\theta \in \Theta$,
\begin{align*}
 | \widehat{l}_{k^* + \ell}(\theta) - \tilde{l}_{k^* + \ell}(\theta) |
 &= | Y_{k^* + \ell}\log f^{k^* + \ell}_{\theta}    - f^{k^* + \ell}_{\theta} - \tilde{Y}_{k^* + \ell}\log \tilde{f}^{k^* + \ell}_{\theta} + \tilde{f}^{k^* + \ell}_{\theta} | \\
 &\leq | Y_{k^* + \ell}\log f^{k^* + \ell}_{\theta} - \tilde{Y}_{k^* + \ell}\log \tilde{f}^{k^* + \ell}_{\theta}| + | f^{k^* + \ell}_{\theta}  - \tilde{f}^{k^* + \ell}_{\theta} |\\
 &\leq    Y_{k^* + \ell} |\log f^{k^* + \ell}_{\theta} - \log \tilde{f}^{k^* + \ell}_{\theta}|
 +  | \log \tilde{f}^{k^* + \ell}_{\theta}| | Y_{k^* + \ell}  - \tilde{Y}_{k^* + \ell} |  + | f^{k^* + \ell}_{\theta}  - \tilde{f}^{k^* + \ell}_{\theta} |\\
.
\end{align*}
By applying the mean value theorem at the function $x \mapsto \log x$ on $[\underline{c}, +\infty [$, we get
$ |\log f^{k^* + \ell}_{\theta} - \log \tilde{f}^{k^* + \ell}_{\theta}| \leq  \frac{1}{\underline{c}} | f^{k^* + \ell}_{\theta}  - \tilde{f}^{k^* + \ell}_{\theta} | $.
Moreover, from the inequality $|\log x | \leq |x-1|, \forall x \geq 1$, we have
$| \log \tilde{f}^{k^* + \ell}_{\theta}| = | \log \frac{ \tilde{f}^{k^* + \ell}_{\theta} }{ \underline{c} }  +  \log \underline{c} |
\leq | \frac{ \tilde{f}^{k^* + \ell}_{\theta} }{ \underline{c} } -1 | + |\log \underline{c} | $.
Hence,
\[ \| \widehat{l}_{k^* + \ell}(\theta) - \tilde{l}_{k^* + \ell}(\theta) \|_{\Theta} \leq C \big( 1 + Y_{k^* + \ell} +   \|\tilde{f}^{k^* + \ell}_{\theta} \|_{\Theta} \big)
          \big( | Y_{k^* + \ell}  - \tilde{Y}_{k^* + \ell} | +  \|f^{k^* + \ell}_{\theta} - \tilde{f}^{k^* + \ell}_{\theta} \|_{\Theta}  \big) ,    \]
 and from the  Hölder's inequality, we get
 \[ \E \| \widehat{l}_{k^* + \ell}(\theta) - \tilde{l}_{k^* + \ell}(\theta) \|^r_{\Theta}
        \leq C \big( \E \big[ 1 + Y_{k^* + \ell} +   \|\tilde{f}^{k^* + \ell}_{\theta} \|_{\Theta} \big]^{\frac{r}{1-r}}  \big)^{1-r}
        \big( \E[ | Y_{k^* + \ell}  - \tilde{Y}_{k^* + \ell} | +  \|f^{k^* + \ell}_{\theta} - \tilde{f}^{k^* + \ell}_{\theta} \|_{\Theta}  ]    \big)^r    .  \]
 From Proposition \ref{prop1}, and arguments of its proof, for all $s>0$, we can find a constant $C>0$ such that $\E Y^s_{k^* + \ell} \leq C$ and
 $ \E \|\tilde{f}^{k^* + \ell}_{\theta} \|^s_{\Theta} \leq C$. Therefore,
  \begin{equation} \label{proof_lem2_Elr}
  \E \| \widehat{l}_{k^* + \ell}(\theta) - \tilde{l}_{k^* + \ell}(\theta) \|^r_{\Theta}
        \leq  C \big( \E | Y_{k^* + \ell}  - \tilde{Y}_{k^* + \ell} | +  \E \|f^{k^* + \ell}_{\theta} - \tilde{f}^{k^* + \ell}_{\theta} \|_{\Theta}    \big)^r    .
  \end{equation}
  According to assumption \textbf{A}$_0 (\Theta)$ and (\ref{lem2_approxYtilde}), we get
 \begin{align*} 
 \E \|f^{k^* + \ell}_{\theta} - \tilde{f}^{k^* + \ell}_{\theta} \|_{\Theta}  &\leq \sum_{j \geq 1} \alpha_j^{(0)} \E | Y_{k^* + \ell -j} - \tilde{Y}_{k^* + \ell -j} |
  \leq \sum_{j = 1}^{\ell/2 - 1} \alpha_j^{(0)} \E | Y_{k^* + \ell -j} - \tilde{Y}_{k^* + \ell -j} |  + C \sum_{j  \geq \ell/2 } \alpha_j^{(0)} \\
  &\leq C \sum_{j = 1}^{\ell/2 - 1} \alpha_j^{(0)} \Big( \inf_{1 \leq p \leq \ell-j } \big\{ (\alpha^{(0)})^{(\ell-j)/p} +  \sum_{ i \geq p} \alpha_i^{(0)} \big\} \Big)
       + C \sum_{j  \geq \ell/2 } \alpha_j^{(0)} \\
 &\leq C \Big(  \inf_{1 \leq p \leq \ell/2 } \big\{ (\alpha^{(0)})^{\ell/(2p)} +  \sum_{ i \geq p} \alpha_i^{(0)} \big\}
       + C \sum_{j  \geq \ell/2 } \alpha_j^{(0)} \Big).
 \end{align*}
 Thus, (\ref{proof_lem2_Elr}) and (\ref{lem2_approxYtilde}) imply
   \begin{align*} 
  \E \| \widehat{l}_{k^* + \ell}(\theta) - \tilde{l}_{k^* + \ell}(\theta) \|^r_{\Theta}
        & \leq   C \Big( \inf_{1 \leq p \leq \ell } \big\{ (\alpha^{(0)})^{\ell/p} +  \sum_{ j \geq p} \alpha_j^{(0)} \big\}
         +   \inf_{1 \leq p \leq \ell/2 } \big\{ (\alpha^{(0)})^{\ell/(2p)} +  \sum_{ j \geq p} \alpha_j^{(0)} \big\}
       + \sum_{j  \geq \ell/2 } \alpha_j^{(0)}      \Big)^r  \\
   &  \leq   C \Big(  \inf_{1 \leq p \leq \ell/2 } \big\{ (\alpha^{(0)})^{\ell/(2p)} +  \sum_{ j \geq p} \alpha_j^{(0)} \big\} + \sum_{j  \geq \ell/2 } \alpha_j^{(0)} \Big)^r \\
   &  \leq   C \Big(   (\alpha^{(0)})^{\ell/(2p_{\ell})} +  \sum_{ j \geq p_{\ell}  } \alpha_j^{(0)}  \Big)^r
    \leq C \Big(   (\alpha^{(0)})^{\ell/(2rp_{\ell})} +  \big( \sum_{ j \geq p_{\ell}  } \alpha_j^{(0)} \big)^r  \Big)
  \end{align*}
 with $p_{\ell} = \ell / \log \ell$ .
 Hence,
 \begin{align}
  \nonumber \sum_{\ell \geq 1} \big( \frac{1}{\ell} \big)^r \E \big[ \| \widehat{l}_{k^* + \ell}(\theta) - \tilde{l}_{k^* + \ell}(\theta) \|^r_{\Theta}  \big]
  \nonumber  & \leq C  \sum_{\ell \geq 1} \big( \frac{1}{\ell} \big)^r  \Big(   (\alpha^{(0)})^{\ell/(2rp_{\ell})} +  \big( \sum_{ j \geq p_{\ell}  } \alpha_j^{(0)} \big)^r  \Big) \\
  \nonumber  & \leq C  \sum_{\ell \geq 1}  \frac{1}{\ell^r}  (\alpha^{(0)})^{\frac{\log \ell}{2r} }
         +  C \sum_{\ell \geq 1}  \frac{1}{\ell^r} \Big(\sum_{ j \geq \ell / \log \ell  } \alpha_j^{(0)} \Big)^r  \\
  \nonumber  & \leq C  \sum_{\ell \geq 1}  \frac{1}{\ell^{r - \frac{\log \alpha^{(0)} }{2r} } }
         +  C \sum_{\ell \geq 1}  \frac{1}{\ell^r} \Big(   \frac{1}{ \big( \frac{\ell}{\log \ell} \big)^{\gamma-1} }     \Big)^r \\
  \label{proof_lem2_lr_bound_2sum}  & \leq C  \sum_{\ell \geq 1}  \frac{1}{\ell^{r - \frac{\log \alpha^{(0)} }{2r} } }
         +  C \sum_{\ell \geq 1}  \frac{ (\log \ell)^{r(\gamma-1)}   }{\ell^{r \gamma}}
  \end{align}
  If $\log \alpha^{(0)} \leq -\frac{1}{2} $, then $ r - \frac{\log \alpha^{(0)} }{2r} > r + \frac{1}{4r}$ and we can choose for instance $r \in [3/4, 1) $ which ensures that
  each of the sum on the right-hand side of (\ref{proof_lem2_lr_bound_2sum}) is finite (recall that, $\gamma > 3/2$ by assumption).
  On the other hand, if $\log \alpha^{(0)} > -\frac{1}{2} $, then $0 < 1 + 2\log \alpha^{(0)} < 1$, and any
  $r \in [\max(\frac{3}{4}, \frac{1 +  \sqrt{1+ 2\log \alpha^{(0)}} }{2}) ~ ,  1)$ ensure that the sums on the right-hand side of (\ref{proof_lem2_lr_bound_2sum}) are finite.
  Thus, one can find $r \in (0,1)$ such that (\ref{proof_lem2_Kounias_cond}) holds ; which achieves the proof of (\ref{proof_theo_consistent_approxLtilde}) and completes
  the proof of the Theorem.
 \Box

    \paragraph{ Proof of Theorem \ref{theo_AN}}  ~
 For any $1\leq i \leq d$, from the Taylor expansion to the function $ \frac{\partial}{\partial \theta_i} \widehat{L}_n(T_{k^*+1, k^* +n},\theta)$,
 there exists $\theta_{n,i}$ between $ \widehat{\theta}(T_{k^*+1, k^* +n})$ and $\theta^*$ such that
 \[ \frac{\partial}{\partial \theta_i} \widehat{L}_n(T_{k^*+1, k^* +n}, \widehat{\theta}(T_{k^*+1, k^* +n})) = \frac{\partial}{\partial \theta_i} \widehat{L}_n(T_{k^*+1, k^* +n},\theta^*)
     + \frac{\partial^2}{\partial\theta \partial\theta_i} \widehat{L}_n(T_{k^*+1, k^* +n},\theta_{n,i})\cdot ( \widehat{\theta}(T_{k^*+1, k^* +n}) - \theta^*) .\]
 Hence,
   \begin{equation} \label{proof_theo_AN_Taylor_Gn}
  n \widehat{G}_n \cdot ( \widehat{\theta}(T_{k^*+1, k^* +n}) - \theta^*)
   = \frac{\partial}{\partial \theta} \widehat{L}_n(T_{k^*+1, k^* +n},\theta^*)
     - \frac{\partial}{\partial \theta} \widehat{L}_n(T_{k^*+1, k^* +n}, \widehat{\theta}(T_{k^*+1, k^* +n})),
  \end{equation}
 with
 \[  \widehat{G}_n  = - \frac{1}{n} \Big( \frac{\partial^2}{\partial\theta \partial\theta_i} \widehat{L}_n(T_{k^*+1, k^* +n},\theta_{n,i}) \Big)_{1\leq i \leq d}  . \]
 Since $\widehat{\theta}(T_{k^*+1, k^* +n}) \limiten \theta^* $ and $\theta^* \in \ring{\Theta}$, for $n$ large enough, $\widehat{\theta}(T_{k^*+1, k^* +n}) \in \ring{\Theta}$
 and $\frac{\partial}{\partial \theta} \widehat{L}_n(T_{k^*+1, k^* +n}, \widehat{\theta}(T_{k^*+1, k^* +n})) = 0$.
 Therefor, (\ref{proof_theo_AN_Taylor_Gn}) gives
    \begin{equation} \label{proof_theo_AN_Taylor_Ln}
  n \widehat{G}_n \cdot ( \widehat{\theta}(T_{k^*+1, k^* +n}) - \theta^*)  = \frac{\partial}{\partial \theta} \widehat{L}_n(T_{k^*+1, k^* +n},\theta^*).
  \end{equation}
  By going along similar lines as in proof of (\ref{proof_theo_consistent_approxLtilde}), we get

     \begin{equation}\label{proof_theo_AN_approxd2L}
   \frac{ 1 }{n} \Big \| \frac{\partial^2}{\partial \theta \partial \theta'} \widehat{L}\big(T_{k^*+1,k^*+n},\theta\big)
  -  \frac{\partial^2}{\partial \theta \partial \theta'} \tilde{L}\big(T_{k^*+1,k^*+n},\theta\big)  \Big\|_{\Theta}  \overset{\texttt{a.s.}}{\underset{n \to \infty}\longrightarrow} 0 ;
  \end{equation}

       \begin{equation}\label{proof_theo_AN_approxdtransposition}
   \frac{ 1 }{n} \Big \| \sum_{t=1}^{n}\frac{1}{\widehat{f}_{\theta}^t}\big(\frac{\partial}{\partial\theta} \widehat{f}_{\theta}^t \big)
                        \big(\frac{\partial}{\partial\theta} \widehat{f}_{\theta}^t \big)'
   - \sum_{t=1}^{n}\frac{1}{\tilde{f}_{\theta}^t}\big(\frac{\partial}{\partial\theta} \tilde{f}_{\theta}^t \big)
                \big(\frac{\partial}{\partial\theta} \tilde{f}_{\theta}^t \big)'  \Big\|_{\Theta}   \overset{\texttt{a.s.}}{\underset{n \to \infty}\longrightarrow} 0 ;
  \end{equation}
    \begin{equation}\label{proof_theo_AN_approxEdL}
  \E \Big( \frac{ 1 }{\sqrt{n}} \Big \| \frac{\partial}{\partial \theta} \widehat{L}\big(T_{k^*+1,k^*+n},\theta\big)
  -  \frac{\partial}{\partial \theta} \tilde{L}\big(T_{k^*+1,k^*+n},\theta\big)  \Big\|_{\Theta}  \Big) \overset{\texttt{a.s.}}{\underset{n \to \infty}\longrightarrow} 0 ;
  \end{equation}
 From Lemma 7.2 and the proof of Theorem 3.2 of \cite{Doukhan2015}, it follows that :
      \begin{equation}\label{proof_theo_AN_approxd2LEell}
   \frac{ 1 }{n} \Big \|  \frac{\partial^2}{\partial \theta \partial \theta'} \tilde{L}\big(T_{k^*+1,k^*+n},\theta\big)
    -  \E\Big(\frac{\partial^2}{\partial \theta \partial \theta'} \tilde{\ell}_0(\theta) \Big) \Big\|_{\Theta}  \overset{\texttt{a.s.}}{\underset{n \to \infty}\longrightarrow} 0 ;
  \end{equation}
       \begin{equation}\label{proof_theo_AN_approxdtranspositionE}
  \Big\| \sum_{t=1}^{n}\frac{1}{\tilde{f}_{\theta}^t}\big(\frac{\partial}{\partial\theta} \tilde{f}_{\theta}^t \big)
                \big(\frac{\partial}{\partial\theta} \tilde{f}_{\theta}^t \big)'
  -    \E \Big( \frac{1}{\tilde{f}_{\theta}^0}\big(\frac{\partial}{\partial\theta} \tilde{f}_{\theta}^0 \big)
                \big(\frac{\partial}{\partial\theta} \tilde{f}_{\theta}^0 \big)' \Big)   \Big\|_{\Theta}   \overset{\texttt{a.s.}}{\underset{n \to \infty}\longrightarrow} 0 ;
  \end{equation}

      \begin{equation}\label{proof_theo_AN_ESigma}
      \E \Big( \frac{1}{\tilde{f}_{\theta^*}^0}\big(\frac{\partial}{\partial\theta} \tilde{f}_{\theta^*}^0 \big)
                \big(\frac{\partial}{\partial\theta} \tilde{f}_{\theta^*}^0 \big)' \Big)
       = -\E\Big( \frac{\partial^2}{\partial \theta \partial \theta'} \tilde{\ell}_0(\theta^*) \Big) = \widetilde{\Sigma};
  \end{equation}
 \begin{equation}\label{proof_theo_conv_dLSigma}
\frac{1}{\sqrt{n}} \frac{\partial}{\partial \theta} \tilde{L}\big(T_{k^*+1,k^*+n},\theta^* \big)
= \frac{1}{\sqrt{n}}\sum_{t=k^* +1 }^{k^* + n} \frac{\partial}{\partial \theta}  \ell_t(\theta^*) \xrightarrow[n\to+\infty]{\mathcal{D}} \mathcal{N} (0,\widetilde{\Sigma}).
\end{equation}
 According to Theorem \ref{theo_consistent}, (\ref{proof_theo_AN_approxd2L}), (\ref{proof_theo_AN_approxdtransposition}), (\ref{proof_theo_AN_approxd2LEell}),
 (\ref{proof_theo_AN_approxdtranspositionE}) and (\ref{proof_theo_AN_ESigma}), we get
 $ \widehat{G}_n  \overset{\texttt{a.s.}}{\underset{n \to \infty}\longrightarrow}  \widetilde{\Sigma}$ and also
 $ \widehat{\Sigma}_n \overset{\texttt{a.s.}}{\underset{n \to \infty}\longrightarrow}  \widetilde{\Sigma} $.
 Hence, for $n$ large enough, $\widehat{G}_n$ is invertible, therefore in addition to (\ref{proof_theo_AN_Taylor_Ln}), (\ref{proof_theo_AN_approxEdL})
 and (\ref{proof_theo_conv_dLSigma}) it holds that
\begin{multline*}
\sqrt{n} ( \widehat{\theta}(T_{k^*+1, k^* +n}) - \theta^*)= \frac{1}{\sqrt{n}} \widehat{G}_n^{-1} \frac{\partial}{\partial \theta} \widehat{L}_n(T_{k^*+1, k^* +n},\theta^*)
     =  \frac{1}{\sqrt{n}} \widetilde{\Sigma}^{-1} \frac{\partial}{\partial \theta} \widetilde{L}_n(T_{k^*+1, k^* +n},\theta^*) + o_P(1) \\
      \xrightarrow[n\to+\infty]{\mathcal{D}} \mathcal{N} (0,\widetilde{\Sigma}^{-1}) .
\end{multline*}
 \Box

  \medskip

 \medskip

 \medskip

 \medskip

 Let $k> m$ and  $T_{1,m}=\{1,\cdots,m\}$,  $T_{\ell,k}=\{\ell,\ell+1,\cdots,k\}$ with $\ell \in \Pi_{m,k}=\{m-v_m,v_m+1,\cdots, k-v_m\}$,
define
   \[ C_{k,\ell}:= \sqrt{m}\, \dfrac{k-\ell}{k}\big\|\Sigma^{-1/2}\cdot\big(\widehat{\theta}(T_{\ell,k}) - \widehat{\theta}(T_{1,m})\big)\big\|, \]
   with $\widehat{\theta}$ defined in (\ref{emv}).
\begin{lem} \label{lem1}
 Under the assumptions of Theorem \ref{theo1},
 \[ \sup_{k>m}  ~\max_{ \ell \in \Pi_{m,k}} ~ \dfrac{1}{b((k-\ell)/m)}\, \big| \widehat{C}_{k,\ell} - C_{k,\ell}  \big|=o_P(1) ~ ~  \text{as} ~ n\rightarrow \infty. \]
 \end{lem}
 \begin{dem}
 For any $m\geq 1$, we have
  \[ \sup_{k>m}  ~\max_{ \ell \in \Pi_{m,k}}  ~ \dfrac{1}{b((k-\ell)/m)}\big| \widehat{C}_{k,\ell} - C_{k,\ell}  \big|
   \leq \dfrac{1}{\inf_{s>0}b(s)} ~
  \sup_{k>m}  ~\max_{ \ell \in \Pi_{m,k}} ~ \big| \widehat{C}_{k,\ell} - C_{k,\ell}  \big| . \]
 Therefore, similar arguments as in the proof of Lemma 7.3 of \cite{Doukhan2015} leads to conclusion.
 \end{dem}
 \medskip

  \medskip

%

  \paragraph{ Proof of Theorem \ref{theo1}}  ~

 Recall that
  \[
 P\{\tau(m)<\infty\}  = P\Big\{ ~ \underset{ m<k \leq [Tm] + 1 } {\mbox{sup}} ~ \underset{\ell \in \Pi_{m,k}} {\mbox{max}} \dfrac{\widehat{C}_{k,\ell}}{b((k-\ell)/m)}>1 \Big\}
   \]
 Hence, it suffices to show that
  \begin{equation}\label{proof_theo1_conv}
 \underset{m<k \leq [Tm] + 1} {\mbox{sup}} ~\underset{ \ell \in \Pi_{m,k}}{\mbox{max}} ~ \dfrac{1}{ b((k-\ell)/m) }
 \widehat{C}_{k,\ell} \xrightarrow [m\to +\infty]{\mathcal{D}}   \underset{1< t \leq T} {\mbox{sup}} ~  \underset{1<s<t}  {\mbox{sup}} ~ \dfrac{\|  W_d(s)-sW_d(1)) \|}{t~b(s)} .
    \end{equation}
 According to Lemma \ref{lem1}, it is enough to show that
  \begin{equation}\label{proof_theo1_conv_wh}
 \underset{m<k \leq [Tm] + 1} {\mbox{sup}} ~\underset{ \ell \in \Pi_{m,k}}{\mbox{max}} ~ \dfrac{1}{ b((k-\ell)/m) }
 C_{k,\ell} \xrightarrow [m\to +\infty]{\mathcal{D}}    \underset{1< t \leq T} {\mbox{sup}} ~  \underset{1<s<t}  {\mbox{sup}} ~ \dfrac{\|  W_d(s)-sW_d(1)) \|}{t~b(s)} .
    \end{equation}

  \medskip
 Let $k>m$ and $\ell \in \Pi_{m,k}$. From the proof of Theorem 4.1 of \cite{Doukhan2015}, it holds that, as $m \rightarrow \infty$
  \[ \Sigma (\widehat{\theta}(T_{1,m}) - \theta_0^*) = \frac{1}{m} \frac{\partial}{\partial \theta} L_m(T_{1m},\theta_0^*) + o_P(\frac{1}{\sqrt{m}}) ~  \text{ and } ~
   \Sigma (\widehat{\theta}(T_{\ell, k}) - \theta_0^*) = \frac{1}{k-\ell} \frac{\partial}{\partial \theta} L_m(T_{\ell,k},\theta_0^*) + o_P(\frac{1}{\sqrt{k-\ell}}) .\]
 Therefore,
  \[ \Sigma (\widehat{\theta}(T_{\ell, k}) - \widehat{\theta}(T_{1, m}) ) = \frac{1}{k-\ell}  \Big( \frac{\partial}{\partial \theta} L_m(T_{\ell,k},\theta_0^*)
           - \frac{k-\ell}{m} \frac{\partial}{\partial \theta} L_m(T_{1,m},\theta_0^*)  \Big) +   o_P( \frac{1}{\sqrt{k-\ell}} + \frac{1}{\sqrt{m}}  ) .\]
 This implies
  \[ C_{k,\ell}  = \frac{ \sqrt{m}}{k} \Sigma^{-1/2} \Big( \frac{\partial}{\partial \theta} L_m(T_{\ell,k},\theta_0^*)
           - \frac{k-\ell}{m} \frac{\partial}{\partial \theta} L_m(T_{1,m},\theta_0^*)  \Big) +   o_P(1) .\]
 Hence,
 \begin{multline*}
  \underset{m<k \leq [Tm] + 1} {\mbox{sup}} ~\underset{ \ell \in \Pi_{m,k}}{\mbox{max}} ~ \dfrac{1}{ b((k-\ell)/m) }
  \Big\| C_{k,\ell} - \frac{ \sqrt{m}}{k} \Sigma^{-1/2} \Big( \frac{\partial}{\partial \theta} L_m(T_{\ell,k},\theta_0^*)
           - \frac{k-\ell}{m} \frac{\partial}{\partial \theta} L_m(T_{1,m},\theta_0^*)  \Big)    \Big\| \\
  \leq  \dfrac{1}{\inf_{s>0}b(s)}  \underset{m<k \leq [Tm] + 1} {\mbox{sup}} ~\underset{ \ell \in \Pi_{m,k}}{\mbox{max}}
  \Big\| C_{k,\ell} - \frac{ \sqrt{m}}{k} \Sigma^{-1/2} \Big( \frac{\partial}{\partial \theta} L_m(T_{\ell,k},\theta_0^*)
           - \frac{k-\ell}{m} \frac{\partial}{\partial \theta} L_m(T_{1,m},\theta_0^*)  \Big)    \Big\| = o_P(1).
 \end{multline*}
 Thus, to complete the proof of the theorem, we will prove that
   \begin{multline}\label{proof_theo1_conv_dL}
 \underset{m<k \leq [Tm] + 1} {\mbox{sup}} ~\underset{ \ell \in \Pi_{m,k}}{\mbox{max}} ~ \dfrac{1}{ b((k-\ell)/m) }
 \frac{ \sqrt{m}}{k} \Big \| \Sigma^{-1/2}  \Big( \frac{\partial}{\partial \theta} L_m(T_{\ell,k},\theta_0^*)
           - \frac{k-\ell}{m} \frac{\partial}{\partial \theta} L_m(T_{1,m},\theta_0^*)  \Big) \Big \| \\
 \xrightarrow [m\to +\infty]{\mathcal{D}}    \underset{1< t \leq T} {\mbox{sup}} ~  \underset{1<s<t}  {\mbox{sup}} ~ \dfrac{\|  W_d(s)-sW_d(1)) \|}{t~b(s)} .
    \end{multline}
 Let $k>m$ and $\ell \in \Pi_{m,k} $. We have
  \[  \dfrac{ \sqrt{m} }{k} \big(  \dfrac{\partial}{\partial\theta}L_m(T_{\ell,k},\theta^*_0) - \dfrac{k-\ell}{m} \dfrac{\partial}{\partial\theta}L_m(T_{1,m},\theta^*_0) \big)
  = - \dfrac{m}{k} \dfrac{1}{\sqrt{m}} \big(  \sum_{i=\ell}^k \dfrac{\partial l_i(\theta^*_0)}{\partial\theta}
  -  \dfrac{k-\ell}{m} \sum_{i=1}^m \dfrac{\partial l_i(\theta^*_0)}{\partial\theta}  \big). \]
 Let us consider the following cases.

 \medskip

 (i) Closed-end procedure. \\
 Let $1<T< \infty$.  Define the set $ S:=\{ (t,s)\in [1,T]\times [1,T] / ~ s<t  \}  $.
 According to \cite{Doukhan2015}, $\big( \dfrac{\partial l_i(\theta^*_0)}{\partial\theta} , \mathcal{F}_i \big)_{i \in \Z}$  is a stationary ergodic  square integrable
 martingale difference sequence with covariance matrix $\Sigma$.
 By the Cramér-Wold device (see \cite{Billingsley1968}), it holds that
     \[  \dfrac{1}{\sqrt{m}}\sum_{i=[ms]}^{[mt]} \dfrac{\partial l_i(\theta^*_0)}{\partial\theta}
     ~ ~ \overset{ \mathcal{D}(S)}{\underset{m\rightarrow \infty} \longrightarrow} ~~  W_{\Sigma}(t-s) \]
 where $\overset{ \mathcal{D}(S)}{\underset{m\rightarrow \infty} \longrightarrow}$ denotes the weak convergence on the Skorohod space $\mathcal{D}(S)$
 and
 $W_{\Sigma}$ is a centered Gaussian process such that $\E \big( W_{\Sigma}(s) , W_{\Sigma}(\tau)'  \big) = \min(s,\tau)\Sigma$.
 Therefore
 \[ \dfrac{1}{\sqrt{m}} \big( \sum_{i=[ms]}^{[mt]} \dfrac{\partial  l_i(\theta^*_0)}{\partial\theta}
 -  \dfrac{[mt]-[ms]}{m}\sum_{i=1}^{m} \dfrac{\partial  l_i(\theta^*_0)}{\partial\theta} \big) \\
   \overset{ \mathcal{D}(S)}{\underset{m\rightarrow \infty} \longrightarrow} ~~  W_{\Sigma}(t-s)-(t-s)B_{\Sigma}(1);
  \]
and
 \[ \dfrac{1}{\sqrt{m}} \Sigma^{-1/2} \big( \sum_{i=[ms]}^{[mt]} \dfrac{\partial  l_i(\theta^*_0)}{\partial\theta}
 -  \dfrac{[mt]-[ms]}{m}\sum_{i=1}^{m} \dfrac{\partial  l_i(\theta^*_0)}{\partial\theta} \big) \\
   \overset{ \mathcal{D}(S)}{\underset{m\rightarrow \infty} \longrightarrow} ~~  W_d(t-s)-(t-s)B_d(1);
  \]
 Hence
    \begin{multline}
  \underset{m<k<[mT]+1} {\mbox{sup}} ~\underset{ \ell \in
\Pi_{m,k}}{\mbox{max}} ~ \dfrac{1}{ b((k-\ell)/m) }
  \dfrac{ \sqrt{m}}{k}  \Big\| \Sigma^{-1/2} \Big( \dfrac{\partial}{\partial\theta}L_m(T_{\ell,k},\theta^*_0) - \dfrac{k-\ell}{n} \dfrac{\partial}{\partial\theta}L_m(T_{1,n},\theta^*_0) \Big) \Big\|    \\
   \xrightarrow [m\to +\infty]{\mathcal{D}}  \underset{1<t<T} {\mbox{sup}} ~  \underset{1<s<t}  {\mbox{sup}} ~ \dfrac{\| W_d(t-s)-(t-s)W_d(1) \|}{t~b(t-s)}
    \overset{\mathcal{D}}{=}  \underset{1<t<T} {\mbox{sup}} ~  \underset{1<s<t}  {\mbox{sup}} ~ \dfrac{\| W_d(s)-s\, W_d(1) \|}{t~b(s)}.
    \end{multline}

 \medskip

 (ii) Open-end procedure. We proceed as in proof of Lemma 6.3 of \cite{Bardet2014}.
      Thus, according to (\ref{proof_theo1_conv_dL}) and (i), it suffices to show that the limit in distribution (as $m, T \rightarrow \infty$) of
      \[ \underset{ k >[Tm] } {\mbox{sup}} ~\underset{ \ell \in \Pi_{m,k}}{\mbox{max}} ~ \dfrac{1}{ b((k-\ell)/m) }
 \frac{ \sqrt{m}}{k} \Big \| \Sigma^{-1/2}  \Big( \frac{\partial}{\partial \theta} L_m(T_{\ell,k},\theta_0^*)
           - \frac{k-\ell}{m} \frac{\partial}{\partial \theta} L_m(T_{1,m},\theta_0^*)  \Big) \Big \|     \]
 exists and is equal to the limit in distribution (as $T \rightarrow \infty$) of
 \[  \underset{ t >T} {\mbox{sup}} ~  \underset{1<s<t}  {\mbox{sup}} ~ \dfrac{\|  W_d(s)-sW_d(1)) \|}{t~b(s)}  .\]
 %
 %
 %
 Let $k > [mT]$. For some $\ell_k \in \Pi_{m,k}$, we have
 \[  \underset{ \ell \in \Pi_{m,k}}{\mbox{max}} ~ \dfrac{1}{ b((k-\ell)/m) } \frac{ \sqrt{m}}{k} \Big \|  \frac{\partial}{\partial \theta} L_m(T_{\ell,k},\theta_0^*) \Big \|
  =  \dfrac{1}{ b((k-\ell_k)/m) } \frac{ \sqrt{m}}{k} \Big \| \sum_{i=\ell_k}^k \dfrac{\partial l_i(\theta^*_0)}{\partial\theta}  \Big \| .\]
 From the Hájek-Rényi-Chow inequality (see Chow (1960) \cite{Chow1960}), we get
 \begin{equation} \label{proof_theo1_Chow}
 \forall x>0, ~  \lim_{T \rightarrow \infty} \limsup_{m \rightarrow \infty } P\Big(  \underset{ k >[Tm] } {\mbox{sup}} ~
 \dfrac{1}{ b((k-\ell_k)/m) } \frac{ \sqrt{m}}{k} \Big \| \sum_{i=\ell_k}^k \dfrac{\partial l_i(\theta^*_0)}{\partial\theta}  \Big \|   > x \Big)=0 .
 \end{equation}
  Moreover, since the function $b(\cdot)$ is non-increasing, we have for any $m,T>1$
   \begin{align}
  \underset{k>mT} {\mbox{sup}} ~\underset{ \ell \in \Pi_{m,k}}{\mbox{max}} ~ \dfrac{1}{ b((k-\ell)/m) } \dfrac{1\sqrt{m}}{k}
  \Big\|  \dfrac{k-\ell}{m}\dfrac{\partial}{\partial\theta}L_m(T_{1,m},\theta^*_0) \Big\|
 \nonumber    &=  \Big\| \dfrac{1}{\sqrt{m}} \sum_{i=1}^m \dfrac{\partial l_i(\theta^*_0)}{\partial\theta} \Big\|
             \times \underset{k>mT} {\mbox{sup}} ~\underset{ \ell \in \Pi_{m,k}}{\mbox{max}} \dfrac{1}{ b((k-\ell)/m) }  \dfrac{k-\ell}{k}  \\
 \nonumber   & =  \Big\| \dfrac{1}{\sqrt{m}} \sum_{i=1}^m \dfrac{\partial l_i(\theta^*_0)}{\partial\theta} \Big\|
             \times \underset{k>mT} {\mbox{sup}} ~ \dfrac{1}{ b((k-v_m)/m) }  \dfrac{k-v_m}{k}  \\
 \nonumber    &=  \dfrac{1}{\underset{ s>0}{\mbox{Inf}}~b(s)} \Big\| \dfrac{1}{\sqrt{m}} \sum_{i=1}^m \dfrac{\partial l_i(\theta^*_0)}{\partial\theta} \Big\|    \\
  \label{proof_theo1_conv_T1m} & \xrightarrow [m\to +\infty]{\mathcal{D}}  \dfrac{1}{\underset{ s>0}{\mbox{Inf}}~b(s)} \|W_{\Sigma}(1) \| ,
     \end{align}
  where the latter convergence holds from the Cramèr-Wold device and the central limit theorem applied to the martingale difference sequence $\big( \dfrac{\partial l_i(\theta^*_0)}{\partial\theta} , \mathcal{F}_i \big)_{i \in \Z}$.
 According to (\ref{proof_theo1_Chow}) and (\ref{proof_theo1_conv_T1m}), it follows that
   \begin{equation}\label{proof_theo1_conv_TlkT1m}
  \underset{k>mT} {\mbox{sup}} ~\underset{ \ell \in \Pi_{m,k}}{\mbox{max}} ~ \dfrac{1}{ b((k-\ell)/m) }
  \dfrac{1\sqrt{n}}{k} \Big\| \Sigma^{-1/2} \Big( \dfrac{\partial}{\partial\theta}L_m(T_{\ell,k},\theta^*_0) - \dfrac{k-\ell}{m} \dfrac{\partial}{\partial\theta}L_m(T_{1,m},\theta^*_0) \Big) \Big\|
 ~ \xrightarrow [m\to +\infty]{\mathcal{D}} ~   \dfrac{1}{\underset{ s>0}{\mbox{Inf}}~b(s)} \|W_d(1) \|.
     \end{equation}
 On the other hand, form the proof of Lemma 6.3 of \cite{Bardet2014}, we get
  \[ \underset{t>T}{\mbox{sup}} ~ \underset{1<s<t}{\mbox{sup}} \dfrac{\| W_{\Sigma}(s) - s W_{\Sigma}(1) \|}{t~b(s)}
  ~ \xrightarrow [m\to +\infty]{\mathcal{D}} ~ \dfrac{1}{\underset{ s>0}{\mbox{Inf}}b(s)} \|W_{\Sigma}(1) \|. \]
 This implies
  \begin{equation}\label{proof_theo1_conv_Wd}
  \underset{t>T}{\mbox{sup}} ~ \underset{1<s<t}{\mbox{sup}} \dfrac{\| W_d(s) - s W_d(1) \|}{t~b(s)}
  ~ \xrightarrow [m\to +\infty]{\mathcal{D}} ~ \dfrac{1}{\underset{ s>0}{\mbox{Inf}}b(s)} \|W_d(1) \|.
     \end{equation}
 (\ref{proof_theo1_conv_TlkT1m}) and (\ref{proof_theo1_conv_Wd}) complete the proof in the case of the open-end procedure.
    \Box

 \paragraph{ Proof of Theorem \ref{theo2} }  ~
 \medskip
 In the sequel, $C$ denotes a positive constant whom value may differ from an inequality to another. \\
  \noindent Denote $k_m= k^* + m^\delta$ for $\delta \in (1/2,1)$.
 For $m$ large enough, we have $m\leq k_m \leq [Tm] + 1$ for both open-end and closed-end procedure; moreover, $v_n << n^\delta$ and $k^* \in \Pi_{m,k_m} $.
 Hence, according to assumption \textbf{B}, we can find a constant $C>0$ such that
   \begin{align}
   \nonumber \underset{ \ell \in \Pi_{m,k_m}}{\max} \frac{\widehat{C}_{k_m,\ell}}{ b((k_m-\ell)/m)}
            &=  \underset{ \ell \in \Pi_{m,k_m}}{\max} \dfrac{1}{ b((k_m-\ell)/m) }
            \sqrt{m} \frac{k_m-\ell}{k_m} \big\| \widehat{\Sigma}_m^{-1/2}   \big(\widehat{\theta}(T_{\ell,k_m}) - \widehat{\theta}(T_{1,m})\big)\big\|  \\
 \nonumber & \geq  \dfrac{1}{ b((k_m-k^*)/m) }
            \sqrt{m} \frac{k_m-k^*}{k_m} \big\| \widehat{\Sigma}_m^{-1/2}   \big(\widehat{\theta}(T_{k^*,k_m}) - \widehat{\theta}(T_{1,m})\big)\big\|   \\
 \nonumber & \geq  C \sqrt{m} ~ \dfrac{m^\delta}{[T^* m] + m^\delta} \big\| \widehat{\Sigma}_m^{-1/2}   \big(\widehat{\theta}(T_{k^*,k_m}) - \widehat{\theta}(T_{1,m})\big)\big\| \\
 \label{proof_theo2_C_sup} & \geq  C ~ m^{\delta - 1/2}  \big\| \widehat{\Sigma}_m^{-1/2}   \big(\widehat{\theta}(T_{k^*,k_m}) - \widehat{\theta}(T_{1,m})\big)\big\| .
   \end{align}
 From \cite{Doukhan2015}, we get $ \widehat{\Sigma}_m^{-1/2} \limitepsm \Sigma^{-1/2}  $ and  $ \widehat{\theta}(T_{1,m}) \limitepsm \theta^*_0  $.
 Moreover, from Theorem \ref{theo_consistent}, $\widehat{\theta}(T_{k^*,k_m}) \limitepsm \theta^*_1 $.
  Thus, since $\Sigma$ is symmetric positive definite,  $\theta^*_0 \neq \theta^*_1$ and $  \delta > 1/2 $, (\ref{proof_theo2_C_sup}) implies
 \[ \underset{ \ell \in \Pi_{m,k_m}}{\max} \frac{\widehat{C}_{k_m,\ell}}{ b((k_m-\ell)/m)}  \limitepsm  \infty  . \]
   \Box

 \bibliographystyle{acm}

\begin{thebibliography}{10}

\bibitem{Ahmad2016}
{\sc Ahmad, A., and Francq, C.}
\newblock Poisson qmle of count time series models.
\newblock {\em Journal of Time Series Analysis 37}, 3 (2016), 291--314.

\bibitem{Bardet2014}
{\sc Bardet, J.-M., and Kengne, W.}
\newblock Monitoring procedure for parameter change in causal time series.
\newblock {\em Journal of Multivariate Analysis 125\/} (2014), 204--221.

\bibitem{Billingsley1968}
{\sc Billingsley, P.}
\newblock Convergence of probability measures.

\bibitem{Chow1960}
{\sc Chow, Y.}
\newblock A martingale inequality and the law of large numbers.
\newblock {\em Proceedings of the American Mathematical Society 11}, 1 (1960),
  107--111.

\bibitem{Chu1996}
{\sc Chu, C.-S.~J., Stinchcombe, M., and White, H.}
\newblock Monitoring structural change.
\newblock {\em Econometrica: Journal of the Econometric Society\/} (1996),
  1045--1065.

\bibitem{Kengne2017}
{\sc Diop, M.~L., and Kengne, W.}
\newblock Testing parameter change in general integer-valued time series.
\newblock {\em Journal of Time Series Analysis 38}, 6 (2017), 880--894.

\bibitem{Kengne2019}
{\sc Diop, M.~L., and Kengne, W.}
\newblock Piecewise autoregression for general integer-valued time series.
\newblock {\em arXiv preprint arXiv:1911.00989\/} (2019).

\bibitem{Doukhan2012}
{\sc Doukhan, P., Fokianos, K., and Tj{\o}stheim, D.}
\newblock On weak dependence conditions for poisson autoregressions.
\newblock {\em Statistics \& Probability Letters 82}, 5 (2012), 942--948.

\bibitem{Doukhan2013corr}
{\sc Doukhan, P., Fokianos, K., Tj{\o}stheim, D., et~al.}
\newblock Correction to "on weak dependence conditions for poisson
  autoregressions"[statist. probab. lett. 82 (2012) 942--948].
\newblock {\em Statistics \& Probability Letters 83}, 8 (2013), 1926--1927.

\bibitem{Doukhan2015}
{\sc Doukhan, P., and Kengne, W.}
\newblock Inference and testing for structural change in general poisson
  autoregressive models.
\newblock {\em Electronic Journal of Statistics 9\/} (2015), 1267--1314.

\bibitem{Doukhan2008}
{\sc Doukhan, P., and Wintenberger, O.}
\newblock Weakly dependent chains with infinite memory.
\newblock {\em Stochastic Processes and their Applications 118}, 11 (2008),
  1997--2013.

\bibitem{Ferland2006}
{\sc Ferland, R., Latour, A., and Oraichi, D.}
\newblock Integer-valued garch process.
\newblock {\em Journal of Time Series Analysis 27}, 6 (2006), 923--942.

\bibitem{Latour2006}
{\sc Ferland, R., Latour, A., and Oraichi, D.}
\newblock Integer-valued garch process.
\newblock {\em Journal of Time Series Analysis 27}, 6 (2006), 923--942.

\bibitem{Franke2012}
{\sc Franke, J., Kirch, C., and Kamgaing, J.~T.}
\newblock Changepoints in times series of counts.
\newblock {\em Journal of Time Series Analysis 33}, 5 (2012), 757--770.

\bibitem{Godambe1960}
{\sc Godambe, V.~P.}
\newblock An optimum property of regular maximum likelihood estimation.
\newblock {\em The Annals of Mathematical Statistics 31}, 4 (1960), 1208--1211.

\bibitem{Gombay2009}
{\sc Gombay, E., and Serban, D.}
\newblock Monitoring parameter change in time series models.
\newblock {\em Journal of Multivariate Analysis 100}, 4 (2009), 715--725.

\bibitem{Horvath2004}
{\sc Horv{\'a}th, L., Hu{\v{s}}kov{\'a}, M., Kokoszka, P., and Steinebach, J.}
\newblock Monitoring changes in linear models.
\newblock {\em Journal of Statistical Planning and Inference 126}, 1 (2004),
  225--251.

\bibitem{Kengne2015}
{\sc Kengne, W.}
\newblock Sequential change-point detection in poisson autoregressive models.
\newblock {\em Journal de la Soci{\'e}t{\'e} Fran{\c{c}}aise de Statistique
  156}, 4 (2015), 98--112.

\bibitem{Kirch2015}
{\sc Kirch, C., and Kamgaing, J.~T.}
\newblock On the use of estimating functions in monitoring time series for
  change points.
\newblock {\em Journal of Statistical Planning and Inference 161\/} (2015),
  25--49.

\bibitem{Kirch2018}
{\sc Kirch, C., Weber, S., et~al.}
\newblock Modified sequential change point procedures based on estimating
  functions.
\newblock {\em Electronic Journal of Statistics 12}, 1 (2018), 1579--1613.

\bibitem{Na2011}
{\sc Na, O., Lee, Y., and Lee, S.}
\newblock Monitoring parameter change in time series models.
\newblock {\em Statistical Methods \& Applications 20}, 2 (2011), 171--199.

\bibitem{Weiss2009}
{\sc Wei{\ss}, C.~H.}
\newblock Modelling time series of counts with overdispersion.
\newblock {\em Statistical Methods and Applications 18}, 4 (2009), 507--519.

\end{thebibliography}

  \end{document}